\documentclass{amsart}
\usepackage{epsfig} 
\usepackage{latexsym} 
\newtheorem{thm}{Theorem}[section]
\newtheorem{cor}[thm]{Corollary}
\newtheorem{lem}[thm]{Lemma}
\newtheorem{prop}[thm]{Proposition}
\theoremstyle{definition}

\theoremstyle{remark}
\newtheorem{rem}[thm]{Remark}
\newtheorem{rems}[thm]{Remarks}
\numberwithin{equation}{section}
\newcommand{\eproof}{\begin{flushright} $\Box$ \end{flushright}}

\newcommand{\ra}{\rightarrow}

\newcommand{\spanm}{\mbox{span}}
\newcommand{\card}[1]{\text{card}\left(#1\right)}
\newcommand{\sing}[1]{\text{sing}\left(#1\right)}

\newcommand{\Q}{\mathbb Q}
\newcommand{\Z}{\mathbb Z}

\renewcommand{\P}{\nabla^{0}}

\usepackage{texdraw}
\usepackage{amsmath}
\def\xplus #1{\bsegment
\setsegscale #1 \move(-10 0) \ravec(20 20) \move(10 0) \rlvec(-9
9) \rmove(-2 2) \ravec(-9 9) \esegment}

\def\xminus #1{\bsegment
\setsegscale #1 \move(-10 0) \rlvec(9 9) \rmove(2 2) \ravec(9 9)
\move(10 0) \ravec(-20 20) \esegment}

\def\xnul #1{\bsegment
\setsegscale #1 \move(-10 0) \ravec(0 20) \move(10 0) \ravec(0 20)
\esegment}

\def\xp #1{\bsegment
\setsegscale #1 \move(-10 0) \ravec(0 20) \move(10 0) \clvec(10
6)(-8 14)(-10 14) \clvec(-14 14)(-14 4)(-10 4) \clvec(-7 4)(-2
6)(0 8) \move(2 10) \clvec(4 12)(8 18)(10 20) \move(9 19) \ravec(1
1) \move(-10 4) \fcir f:0 r:0.6 \move(-10 14) \fcir f:0 r:0.6
\esegment}

\def\xpminus #1{\bsegment
\setsegscale #1 \move(-10 0) \clvec(8 10)(8 10)(-10 20) \move(10
0) \clvec(-8 10)(-8 10)(10 20) \move(0 5.7) \fcir f:0 r:0.6
\move(0 14.3) \fcir f:0 r:0.6 \move(-8 19) \ravec(-2 1) \move(8
19) \ravec(2 1) \esegment}

\def\xpplus #1{\bsegment
\setsegscale #1 \move(-10 0) \ravec(0 20) \move(10 0) \clvec(10
6)(-8 14)(-9 14) \rmove(-2 0) \clvec(-14 14)(-14 4)(-10 4)
\clvec(-7 4)(-2 6)(0 8) \clvec(4 12)(8 18)(10 20) \move(9 19)
\ravec(1 1) \move(-10 4) \fcir f:0 r:0.6 \move(1.1 9.1) \fcir f:0
r:0.6

\esegment}

\def\xdot #1{\bsegment
\setsegscale #1 \move(-10 0) \ravec(20 20) \move(10 0) \ravec(-20
20) \move(0 10) \fcir f:0 r:0.8 \esegment}

\def\vparent #1{\bsegment
\setsegscale #1 \clvec(0 0)(-3 10)(0 20) \esegment}

\def\hparent #1{\bsegment
\setsegscale #1 \clvec(0 0)(3 10)(0 20) \esegment}

\def\xdota #1{\bsegment
\setsegscale #1 \xdot{#1} \clvec(0 5)(-3 3)(-5 5) \clvec(-7 7)(-7
9)(-7 10) \clvec(-7 11)(-7 13)(-6 14) \move(-4 16) \clvec(-3 17)(0
15)(0 19) \move(0 19)\ravec(0 1) \move(-5 5) \fcir f:0 r:0.8
\esegment}

\def\xdotc #1{\bsegment
\setsegscale #1 \xdot{#1} \clvec(0 5)(3 3)(5 5) \clvec(7 7)(7 9)(7
10) \clvec(7 11)(7 13)(6 14) \move(4 16) \clvec(3 17)(0 15)(0 19)
\move(0 19)\ravec(0 1) \move(5 5) \fcir f:0 r:0.8 \esegment}

\def\xdotb #1{\bsegment
\setsegscale #1 \move(-10 0) \rlvec(4 4) \move(-4 6) \ravec(14 14)
\move(10 0) \ravec(-20 20) \move(0 10) \fcir f:0 r:0.8 \move(0 0)
\clvec(0 5)(-3 3)(-5 5) \clvec(-7 7)(-7 9)(-7 10) \clvec(-7 11)(-7
13)(-5 15) \clvec(-3 17)(0 15)(0 19) \move(0 19)\ravec(0 1)
\move(-5 15) \fcir f:0 r:0.8 \esegment}

\def\xdotd #1{\bsegment
\setsegscale #1 \move(10 0) \rlvec(-4 4) \move(4 6) \ravec(-14 14)
\move(-10 0) \ravec(20 20) \move(0 10) \fcir f:0 r:0.8 \move(0 0)
\clvec(0 5)(3 3)(5 5) \clvec(7 7)(7 9)(7 10) \clvec(7 11)(7 13)(5
15) \clvec(3 17)(0 15)(0 19) \move(0 19)\ravec(0 1) \move(5 15)
\fcir f:0 r:0.8 \esegment}

\def\xdote #1{\bsegment
\setsegscale #1 \move(-10 0) \rlvec(6.985504245 11.142507074)
\move(-1.985504245 13.357492926) \avec(2 20) \move(-2 0)\ravec(12
20) \move(10 0)\ravec(-20 20) \move(2.5 7.5) \fcir f:0 r:0.8
\esegment}

\def\xdotg #1{\bsegment
\setsegscale #1 \move(10 0) \rlvec(-6.985504245 11.142507074)
\move(1.985504245 13.357492926) \avec(-2 20) \move(2 0)\ravec(-12
20) \move(-10 0)\ravec(20 20) \move(-2.5 7.5) \fcir f:0 r:0.8
\esegment}

\def\xdoth #1{\bsegment
\setsegscale #1 \move(-10 0) \avec(2 20) \move(-2.5 12.5) \fcir
f:0 r:0.8 \move(-2 0)\avec(10 20) \move(10 0)\lvec(3.5 6.5)
\rmove(-2 2) \avec(-10 20) \esegment}

\def\xdotf #1{\bsegment
\setsegscale #1 \move(10 0) \avec(-2 20) \move(2.5 12.5) \fcir f:0
r:0.8 \move(2 0)\avec(-10 20) \move(-10 0)\lvec(-3.5 6.5) \rmove(2
2) \avec(10 20) \esegment}

\def\xminusloop #1{\bsegment
\setsegscale #1 \move(-10 0)\rlvec(1 1) \clvec(-5 5)(-5 15)(-10
10) \clvec(-15 5)(-5 5)(-1 9) \move(-6.3 6.3) \fcir f:0 r:0.8
\move(1 11) \ravec(9 9) \move(10 0) \ravec(-20 20) \esegment}

\def\xnulcirc #1{\bsegment
\setsegscale #1 \move(-10 0) \ravec(0 20) \move(10 0) \ravec(0 20)
\move(-10 9) \larc r:6 sd:270 ed:75 \larc r:6 sd:105 ed:270
\move(-10 3) \fcir f:0 r:0.8 \move(-4 9) \ravec(0 1) \esegment}

\def\xdotaa #1{\bsegment
\setsegscale #1 \xdot{#1} \clvec(0 5)(-3 3)(-5 5) \clvec(-7 7)(-7
9)(-7 10) \clvec(-7 11)(-7 13)(-5 15) \clvec(-3 17)(0 15)(0 19)
\move(0 19)\ravec(0 1) \move(-5 5) \fcir f:0 r:0.8 \esegment}

\def\xdotcc #1{\bsegment
\setsegscale #1 \xdot{#1} \clvec(0 5)(3 3)(5 5) \clvec(7 7)(7 9)(7
10) \clvec(7 11)(7 13)(5 15) \clvec(3 17)(0 15)(0 19) \move(0
19)\ravec(0 1) \move(5 5) \fcir f:0 r:0.8 \esegment}

\def\xdotbb #1{\bsegment
\setsegscale #1 \move(-10 0) \rlvec(5 5) \move(-5 5) \ravec(15 15)
\move(10 0) \ravec(-20 20) \move(0 10) \fcir f:0 r:0.8 \move(0 0)
\clvec(0 5)(-3 3)(-5 5) \clvec(-7 7)(-7 9)(-7 10) \clvec(-7 11)(-7
13)(-5 15) \clvec(-3 17)(0 15)(0 19) \move(0 19)\ravec(0 1)
\move(-5 15) \fcir f:0 r:0.8 \esegment}

\def\xdotdd #1{\bsegment
\setsegscale #1 \move(10 0) \rlvec(-5 5) \move(5 5) \ravec(-15 15)
\move(-10 0) \ravec(20 20) \move(0 10) \fcir f:0 r:0.8 \move(0 0)
\clvec(0 5)(3 3)(5 5) \clvec(7 7)(7 9)(7 10) \clvec(7 11)(7 13)(5
15) \clvec(3 17)(0 15)(0 19) \move(0 19)\ravec(0 1) \move(5 15)
\fcir f:0 r:0.8 \esegment}

\def\xdotee #1{\bsegment
\setsegscale #1 \move(-10 0) \avec(2 20) \move(-2 0)\ravec(12 20)
\move(10 0)\ravec(-20 20) \move(2.5 7.5) \fcir f:0 r:0.8
\esegment}

\def\xdotgg #1{\bsegment
\setsegscale #1 \move(10 0) \avec(-2 20) \move(2 0)\ravec(-12 20)
\move(-10 0)\ravec(20 20) \move(-2.5 7.5) \fcir f:0 r:0.8
\esegment}

\def\xdothh #1{\bsegment
\setsegscale #1 \move(-10 0) \avec(2 20) \move(-2.5 12.5) \fcir
f:0 r:0.8 \move(-2 0)\avec(10 20) \move(10 0) \avec(-10 20)
\esegment}

\def\xdotff #1{\bsegment
\setsegscale #1 \move(10 0) \avec(-2 20) \move(2.5 12.5) \fcir f:0
r:0.8 \move(2 0)\avec(-10 20) \move(-10 0) \avec(10 20) \esegment}

\def\xdotccc #1{\bsegment
\setsegscale #1 \xdot{#1} \move(5 5) \clvec(7 7)(7 9)(7 10)
\clvec(7 11)(7 13)(6 14) \move(4 16) \clvec(3 17)(0 15)(0 19)
\move(0 19)\ravec(0 1) \move(5 5) \fcir f:0 r:0.8 \move(0 20)
\clvec(0 30)(-18 25)(-18 10) \clvec(-18 -10)(-13 -3)(-10 0)

\move(0 -17) \rlvec(4 4) \move(10 -17) \rlvec(-5 5) \move(5 -12)
\clvec(0 -7)(0 0)(5 5) \move(6 -11) \clvec(8 -9)(15 -5)(10 0)
\esegment}

\def\xdotbbb #1{\bsegment
\setsegscale #1 \move(-10 0) \rlvec(4 4) \move(-4 6) \ravec(14 14)
\move(10 0) \ravec(-20 20) \move(0 10) \fcir f:0 r:0.8 \move(-5 5)
\clvec(-7 7)(-7 9)(-7 10) \clvec(-7 11)(-7 13)(-5 15) \clvec(-3
17)(0 15)(0 19) \move(0 19)\ravec(0 1) \move(-5 15) \fcir f:0
r:0.8

\move(0 20) \clvec(0 30)(-18 25)(-18 10) \clvec(-18 -10)(-13
-3)(-10 0) \move(0 -17) \rlvec(4 4) \move(10 -17) \rlvec(-5 5)
\move(5 -12) \clvec(2 -9)(-3 3)(-5 5) \move(6 -11) \clvec(8 -9)(15
-5)(10 0) \esegment}

\def\xdotddd #1{\bsegment
\setsegscale #1 \move(10 0) \rlvec(-4 4) \move(4 6) \ravec(-14 14)
\move(-10 0) \ravec(20 20) \move(0 10) \fcir f:0 r:0.8 \move(5 5)
\clvec(7 7)(7 9)(7 10) \clvec(7 11)(7 13)(5 15) \clvec(3 17)(0
15)(0 19) \move(0 19)\ravec(0 1) \move(5 15) \fcir f:0 r:0.8
\move(0 20) \clvec(0 30)(-18 25)(-18 10) \clvec(-18 -10)(-13
-3)(-10 0)

\move(0 -17) \rlvec(4 4) \move(10 -17) \rlvec(-5 5) \move(5 -12)
\clvec(0 -7)(0 0)(5 5) \move(6 -11) \clvec(8 -9)(15 -5)(10 0)
\esegment}

\def\xdoteee #1{\bsegment
\setsegscale #1 \move(-10 0) \rlvec(6.985504245 11.142507074)
\move(-1.985504245 13.357492926) \avec(2 20) \move(2 20) \clvec(8
26)(-10 30)(-14 26) \move(-16 24) \clvec(-20 20)(-18 -8)(-10 0)
\move(2.5 7.5)\avec(10 20) \move(2.5 7.5)\avec(-20 30) \move(2.5
7.5) \fcir f:0 r:0.8

\move(10 -12) \rlvec(-6 6) \move(-2 -12) \rlvec(5 5) \move(5 -5)
\clvec(9 -1)(5.5 4.5)(2.5 7.5) \move(4 -6) \clvec(0 -2)(-0.5
4.5)(2.5 7.5) \esegment}

\def\xdotfff #1{\bsegment
\setsegscale #1 \move(2.5 12.5) \avec(-2 20) \move(2.5 12.5) \fcir
f:0 r:0.8 \move(10 -13.333333333)\avec(-16 30)

\move(2 -13.333333333) \lvec(5.1 -8.1666666666) \move(6.9
-5.166666666) \clvec(9.9 -0.166666666)(5.5 7.5)(2.5 12.5)
\move(-10 0)\lvec(-3.5 6.5) \rmove(2 2) \avec(10 20)

\move(-2 20) \clvec(-5 25)(-10 28)(-12 26) \move(-14 24)
\clvec(-18 20)(-18 -8)(-10 0) \esegment}

\def\xdotggg #1{\bsegment
\setsegscale #1 \move(10 0) \move(1.985504245 13.357492926)
\avec(-2 20) \move(10 -13.333333333)\avec(-16 30) \move(2
-13.333333333) \lvec(5.1 -8.1666666666) \move(6.9
-5.166666666)\clvec(9.9 -0.166666666)(6.1 6.5)(3.1 11.5) \move(-10
0)\lvec(-3.5 6.5) \move(-10 0)\ravec(20 20) \move(-2.5 7.5) \fcir
f:0 r:0.8 \move(-2 20) \clvec(-5 25)(-10 28)(-12 26) \move(-14 24)
\clvec(-18 20)(-18 -8)(-10 0) \esegment}

\def\xdothhh #1{\bsegment
\setsegscale #1 \move(-10 0) \avec(2 20) \move(-2.5 12.5) \fcir
f:0 r:0.8 \move(2.5 7.5)\avec(10 20) \move(3.5 6.5) \rmove(-2 2)
\avec(-20 30) \move(2 20) \clvec(8 26)(-10 30)(-14 26) \move(-16
24) \clvec(-20 20)(-18 -8)(-10 0)

\move(10 -12) \rlvec(-6 6) \move(-2 -12) \rlvec(5 5) \move(5 -5)
\clvec(9 -1)(5.5 4.5)(3.5 6.5) \move(4 -6) \clvec(0 -2)(-0.5
4.5)(2.5 7.5) \esegment}

\def\xnulcircv #1{\bsegment
\setsegscale #1 \move(-10 0) \ravec(0 20) \move(10 0) \ravec(0 20)
\move(-10 9) \larc r:6 sd:270 ed:90 \larc r:6 sd:90 ed:270
\move(-10 3) \fcir f:0 r:0.8 \move(-4 9) \ravec(0 1) \esegment}

\def\xnulcirch #1{\bsegment
\setsegscale #1 \move(-10 0) \ravec(0 20) \move(10 0) \ravec(0 20)
\move(10 9) \larc r:6 sd:270 ed:90 \larc r:6 sd:90 ed:270 \move(10
3) \fcir f:0 r:0.8 \move(16 9) \ravec(0 1) \esegment}

\begin{document}

\title[Higher Skein Modules]{Higher Skein Modules}%
\author{J{\o}rgen Ellegaard Andersen}%
\address{Department of Mathematics, Aarhus University, DK-8000 Aarhus C}%
\email{andersen@math.au.dk}%
\author{Vladimir Turaev}%
\address{Institut de Recherche Math\'{e}matique Avanc\'{e}e, Universit\'{e}
Louis Pastuer-CNRS,
7 rue Ren\'{e} Descartes, 67084 Strasbourg Cedex, France.}%
\email{turaev@math.u-strasbg.fr}%

\thanks{This work was supported
by the TMR network 'Algebraic Lie Representaions ' of the European
Commission, EC Contract No ERB FMRX CT97-0100 and by MaPhySto --
centre for Mathematical Physics and Stochastics, funded by The
Danish National Research Foundation.}

\thanks{Joint Aarhus University and MaPhySto preprint.}

\date{November 23, 1998}%
\begin{abstract}
We introduce higher skein modules of links generalizing the Conway
skein module.
  We show that these modules are closely connected to
 the HOMFLY polynomial.
\end{abstract}
\maketitle
\section{Introduction}\label{intro}

\begin{figure}[htbp]
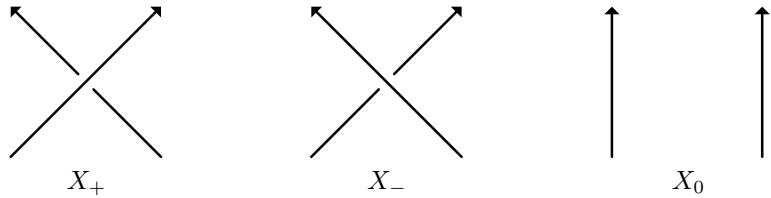

  \centertexdraw{\drawdim mm
\linewd 0.3 \arrowheadtype t:F \arrowheadsize l:1 w:2 \xplus{1}
\move(40 0) \xminus{1} \rmove(40 0) \xnul{1} \move(0 -2) \textref
h:C v:T \htext{$X_{+}$} \move(40 -2) \textref h:C v:T
\htext{$X_{-}$} \move(80 -2) \textref h:C v:T \htext{$X_{0}$}
\move(0 -10) }
  \caption{  $X_+,X_-$, and $X_0$.}\label{fot}
\end{figure}

The notion of a skein module  of links
arises naturally from the study of  link polynomials.
For instance, the  one-variable Conway polynomial   $\nabla $
of   links in the 3-sphere $S^3$ satisfies the fundamental
skein relation $\nabla (X_+)=\nabla (X_-)+h\nabla (X_0)$
 where $h$ is the variable and
$X_+,X_-,X_0$ are any three oriented links coinciding outside a
3-ball and looking as in Figure \ref {fot} inside this ball. This
suggests to consider the $\Z[h]$-module generated by the isotopy
classes of oriented links in $S^3$ modulo the relations
$X_+-X_--hX_0=0$ corresponding to all triples $(X_+,X_-,X_0)$ as
above. This is the Conway skein module of $S^3$. Applying similar
definitions to links in an oriented 3-manifold $M$, we obtain the
skein module of $M$, cf. \cite{Pr} and \cite{Tu}.

\begin{figure}[htbp]
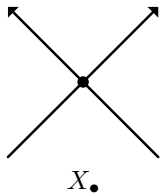

  \centertexdraw{\drawdim mm
\linewd 0.3 \arrowheadtype t:F \arrowheadsize l:1 w:2 \xdot{1}
\move(0 -2) \textref h:C v:T \htext{$X_{\bullet}$} \move(0 -10) }
  \caption{A double point $X_\bullet$.}\label{doublept}
\end{figure}

In this paper we introduce \lq\lq higher" versions of the
 Conway skein module. Our approach is    inspired by the theory of Vassiliev
link invariants.
 In that theory one considers singular links, i.e.,
links with  double points as in Figure \ref {doublept}. Each double point
  $X_\bullet$   is
resolved into the formal difference $X_+-X_-$ of the positive and negative
crossings. In this way each   singular link with $n$ double points is
resolved
into a formal linear combination of $2^n$ (non-singular) links. This yields
the
Vassiliev filtration
  $V=V_0\supset V_1\supset V_2\supset \ldots$
  where $V$ is the abelian group freely generated by the isotopy classes of
oriented
links   and $V_n$ is its subgroup generated by the resolutions of
   singular links with $n$ double points. The  consecutive quotients
   $V_n/V_{n+1}$ are among   the main objects of the theory of
Vassiliev  invariants.

\begin{figure}[htbp]
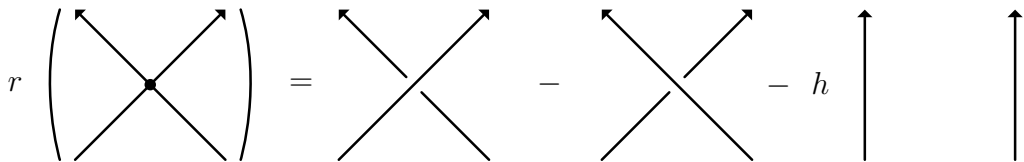

  \centertexdraw{\drawdim mm
\linewd 0.3 \arrowheadtype t:F \arrowheadsize l:1 w:2 \move(-12 0)
\vparent{1} \move(0 0)\xdot{1}\rmove(12 0) \hparent{1}\rmove(23
0)\xplus{1}\rmove(35 0)\xminus{1}\rmove(35 0)\xnul{1} \move(-17
10)\textref h:R v:C \htext{\Large $r$} \move(20 10)\textref h:C
v:C \htext{\Large $=$} \move(53 10)\textref h:C v:C \htext{\Large
$-$} \move(86 10)\textref h:C v:C \htext{\Large $-\ \ h$} }
  \caption{The   resolution $r$ of a double point.}\label{qres}
\end{figure}

We consider here a deformation of the Vassiliev filtration. The
idea is to resolve the double points via the formula in Figure
\ref{qres}. More precisely, for
 an oriented 3-manifold $M$ denote by
 $A(M)$   the free $\Z[h]$-module generated by the isotopy
classes of oriented (non-empty) links in  $M$. Here  $\Z[h]$ is
  the ring of polynomials in $h$ with integer coefficients.
By a singular link in $M$, we
mean an immersion of a finite system of oriented circles in
$M$ with only double transversal intersections. Using the
formula
$$r(  X_\bullet  )= X_+-X_- -hX_0,$$
  we resolve each singular link $L$
with $n$ double points into a formal sum $r(L)\in A(M)$ of $3^n$
terms. Denote by $A_n$   the $\Z[h]$-submodule of $A(M)$ generated by $r(L)$
where $L$ runs over all singular links with $n$ double points.
Clearly, $A(M)=A_0\supset A_1\supset A_2\supset \ldots$.
 The quotient $A_0/A_1$ is the Conway skein module of $M$.
We call the $\Z[h]$-modules $A_n/A_{n+1}$ with $n=1,2,\ldots$ the higher
Conway
skein
modules of $M$. Our aim is to compute them  (at least partially)
in the case $M=S^3$. In the sequel, we restrict ourselves to links in $S^3$
unless explicitly stated to the contrary. Set $A=A(S^3)$.

Recall first the structure of $A_0/A_1=A/A_1$. The Conway
polynomial $\nabla$ (normalised so that its value on an unknot is
$1$) defines a $\Z[h]$-linear epimorphism $A/A_1 \to \Z[h]$. Its
kernel
 is  a free abelian group freely generated by the
  classes in $A/A_1$ of the trivial links with $\geq 2$ components.
A standard argument (reproduced below) shows that
 multiplication by $h$ annihilates these classes. Therefore
$A/A_1=h{\text {-torsion}} \oplus \Z[h]$. (By the $h$-torsion of a
$\Z[h]$-module $M$ we mean the set $\{a\in M\,\vert\, ha=0\}$).

\begin{figure}[htbp]
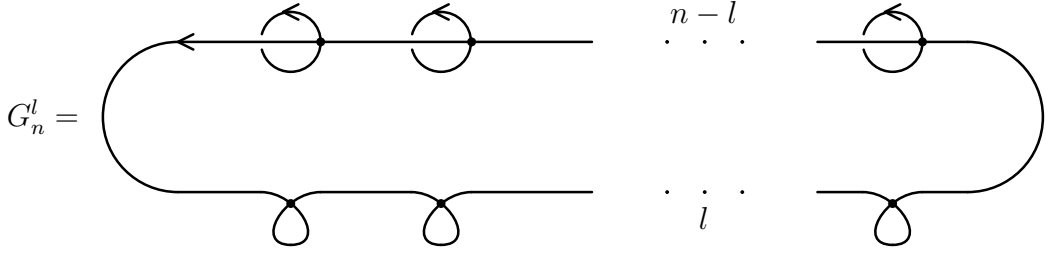

  \centertexdraw{
\drawdim mm \linewd 0.3 \arrowheadtype t:V \arrowheadsize l:2 w:2
\move(0 0) \larc r:10 sd:90 ed:270 \move(105 10) \rlvec(-20 0)
\rmove(-30 0) \ravec(-55 0) \move(105 0) \larc r:10 sd:270 ed:90
\move(115 0) \move(65 10) \fcir f:0 r:0.25 \move(70 10) \fcir f:0
r:0.25 \move(75 10) \fcir f:0 r:0.25 \move(65 -10) \fcir f:0
r:0.25 \move(70 -10) \fcir f:0 r:0.25 \move(75 -10) \fcir f:0
r:0.25

\move(15 10) \larc r:4 sd:0 ed:165 \larc r:4 sd:195 ed:0 \move(19
10) \fcir f:0 r:0.6 \move(15 14) \ravec(-1 0)

\move(35 10) \larc r:4 sd:0 ed:165 \larc r:4 sd:195 ed:0 \move(39
10) \fcir f:0 r:0.6 \move(35 14) \ravec(-1 0)

\move(95 10) \larc r:4 sd:0 ed:165 \larc r:4 sd:195 ed:0 \move(99
10) \fcir f:0 r:0.6 \move(95 14) \ravec(-1 0) \move(0 -10)
\rlvec(11 0)

\move(11 -10) \clvec(16 -10)(20 -17)(15 -17) \move(19 -10)
\clvec(14 -10)(10 -17)(15 -17) \move(15 -11.5) \fcir f:0 r:0.6

\move(19 -10) \rlvec(12 0) \move(31 -10) \clvec(36 -10)(40 -17)(35
-17) \move(39 -10) \clvec(34 -10)(30 -17)(35 -17) \move(35 -11.5)
\fcir f:0 r:0.6

\move(91 -10) \clvec(96 -10)(100 -17)(95 -17) \move(99 -10)
\clvec(94 -10)(90 -17)(95 -17) \move(95 -11.5) \fcir f:0 r:0.6
\move(99 -10) \rlvec(6 0) \move(85 -10) \rlvec(6 0) \move(39 -10)
\rlvec(16 0) \move(-13 0)\textref h:R v:C \htext{\Large
$G_{n}^l=$} \move(70 12)\textref h:C v:B \htext{\Large ${n-l}$}
\move(70 -12)\textref h:C v:T \htext{\Large ${l}$} \move(0 17)
\move(0 -17) \move(-20 0) }
  \caption{The singular link $G^l_{n}$.}\label{Gln}
\end{figure}

It turns out that similar results hold for the higher skein modules.
To state our theorems we introduce a two-parameter family of
singular links $G^l_n$ with $n$ double points where $n\geq 0$ and
$l=0,1,\ldots,n$, see Figure \ref{Gln}. The singular link $G^l_n$
is formed by $n-l$ embedded circles crossing one immersed circle
with $l$ curls. According to our definitions,
  $G^l_n$ represents an element $ r(G^l_n) \in
A_{n}/A_{n+1}$.

Our main theorem on the structure of $A_n/A_{n+1}$ says the
following.

\begin{thm}\label{Mainth}
 For  each $n\geq 0$,
\[A_n/A_{n+1} = h{\text {-torsion}}\, \oplus \,\bigoplus_{l=0}^n \Z[h] \,
r(G^l_n) .\]
\end{thm}

This theorem may be reformulated by saying that the quotient of
the module $A_n/A_{n+1}$ by its $h$-torsion is a free
$\Z[h]$-module of rank $n+1$ with free generators $r(G^l_n) \mod
A_{n+1}$ where $l=0,1,\ldots,n$.
The $h$-torsion of $A_n/A_{n+1}$ is in general quite big. Indeed
setting $h=0$ we obtain a projection $A \to V$ mapping  the filtration
  $A=A_0\supset  A_1\supset
\ldots$ onto the Vassiliev filtration
  $V=V_0\supset V_1\supset   \ldots$.
   This  induces a $\Q$-linear epimorphism
  from $\Q\otimes (A_n/A_{n+1})$ onto the vector space $ \Q\otimes
(V_n/V_{n+1})$
 isomorphic to the vector space of chord
  diagrams with $n$ chords modulo the 4T- and 1T-relations.
It follows from Theorem \ref {Mainth} that
  $\Q\otimes (h-{\text {torsion \,of\,}}  A_n/A_{n+1})$ lies in
the latter vector
space of chord diagrams as a subspace of codimension $n+1$.

  In this paper we shall be mainly interested in the
  $\Z[h]$-free part of $A_n/A_{n+1}$.
Tensoring with the ring of Laurent polynomials $\Lambda=
\Z[h,h^{-1}]$ we obtain the following.

\begin{cor}\label{cor1}
For each $n\geq 0$,
\[  \Lambda \otimes_{\Z[h]} (A_n/A_{n+1})=
\bigoplus_{l=0}^n \Lambda \,r(G^l_n) \,\,\,\,and \,\,\,\, \Lambda
\otimes_{\Z[h]} (A /A_{n})= \bigoplus_{\stackrel{{\scriptstyle
l,m\geq 0}} { l+ m<n}} \Lambda \,r(G^l_{l+m}).
\]
\end{cor}

\begin{cor}\label{linkexp}
There are unique  $\Z[h ]$-linear homomorphisms $\nabla_{l,m}:A\to \Lambda$
numerated by pairs of non-negative
integers $(l,m)$ such that for any $a\in A$,
\[a = \sum_{l,m}  \nabla_{l,m} (a)\,    r(G^l_{l+m}) \in
\varprojlim_n \Lambda \otimes_{\Z[h]} (A/A_n).\]

\end{cor}

Applying this to any oriented link $L$ we obtain an expansion
\[L = \sum_{  l, m}\nabla_{l,m} (L)\,   r(G^l_{l+m}) \in \varprojlim_n
 \Lambda \otimes_{\Z[h]} (A/A_n).\]

By definition of $\nabla_{l,m}$, we have
$$\nabla_{l,m}(A_{l+m+1})=0\,\,\,\,\,\,{\text {and}}\,\,\,\,\,\,
\nabla_{l,m}( r(G^{l'}_{m'})) =\delta_l^{l'} \delta_{l+m}^{m'}$$
where $\delta$ is the Kronecker delta. It is easy to check that
$h^{l+m} \nabla_{l,m}(A) \subset \Z[h]$, for all $l,m$. In
particular, $\nabla_{0,0}$ annihilates $A_1$ and maps the trivial
knot $G^0_0$ into 1. Hence $\nabla_{0,0}=\nabla$ is the Conway
polynomial. The following theorem computes $\nabla_{l,m}$ from
$\nabla_{l,0}$.

\begin{thm}\label{deri}
For any $l\geq 0,m\geq 1$ and any oriented link $L$,
$$\nabla_{l,m}(L)= \frac {h^{-l}}{m!}
(h^{l}\nabla_{l,0}(L))^{(m)}$$ where $f^{(m)}$ is the $m$-th
derivative of a
 Laurent polynomial $f\in \Lambda$. In particular,
$$\nabla_{0,m}(L)= (m!)^{-1} (\nabla (L))^{(m)}.$$

\end{thm}

Theorem \ref{deri} and the inclusion $h^l\nabla_{l,0}(L)\in\Z[h]$
imply the following.

\begin{cor}\label{Poly}
For all links $L$ and all $l,m$ we have that
$h^l\nabla_{l,m}(L)\in\Z[h].$
\end{cor}

It turns out that the sequence of link polynomials $\nabla_{l,0}$ with
$l=0,1,\ldots$
is equivalent to the Homfly
polynomial.  By the Homfly polynomial, $\tilde \nabla$,
we shall mean the (unique)
 mapping  from the set of isotopy classes of oriented links in $S^3$ into
the ring of Laurent polynomials $\Z[x,x^{-1},h,h^{-1}]$
which is uniquely characterised by the following two properties:

(i) the value of $\tilde \nabla$ on an unknot is equal to 1;

(ii) for  any three oriented links   $X_+,X_-,X_0$
coinciding outside a 3-ball
and looking
as in Figure \ref {fot} inside this ball, we have that
\begin{equation}
  x\tilde\nabla (X_+)-x^{-1}\tilde \nabla (X_-)=h\tilde \nabla (X_0).\label{fil1}
 \end{equation}
Clearly, the Conway polynomial $\nabla$ is obtained from $\tilde \nabla$
by the substitution $x=1$.

\begin{thm}\label{thhomfly}
For any oriented link $L$, the formal power series
 $$P(L)(h,u) = \sum_{l\geq 0} (-h)^l\nabla_{l,0}(L)(h)\,u^l$$
is a reparametrisation of $ \tilde
\nabla(L)$.

\end{thm}

The precise form of the reparametrisation   in this theorem is a
little
technical (of course it does not depend on the choice of $L$).
We shall give a   detailed statement in Section 6.
Theorems \ref {deri}
and \ref {thhomfly}
imply that the polynomials $\{\nabla_{l,m} \}_{l,m}$ determine
and are determined by the Homfly polynomial.

 \rems \label{remIntrod}

 1. Applying the definition of the
resolution $r$ inductively to all $l$ curls of $G^l_{l+m}$, we
obtain $$r(G^l_{l+m})= -hu \,r(G^{l-1}_{l+m-1})= h^2 u^2
r(G^{l-2}_{l+m-2})=\ldots=(-h)^l u^l r(G^{0}_{m})$$ where the
variable $u$ acts on $A$ as the disjoint union with an unknot.
(Note that multiplication by $hu$ maps each $A_n$ into $A_{n+1}$).
This gives for each oriented link $L$, an expansion
\[L = \sum_{  l, m}\nabla^0_{l,m} (L)\,   u^l r(G^0_{m}) \in \varprojlim_n
 \Lambda \otimes_{\Z[h]} (A/A_n) \] where
$$\nabla^0_{l,m}(L) =(-h)^{l}\nabla_{l,m}(L) \in \Z[h].$$

2.  It is instructive to set $h=\pm 1$ in our constructions.
To this end, consider the abelian group $V$
freely generated by the isotopy classes of
oriented links. For $\varepsilon= \pm 1$,
 consider the   homomorphism  $ A\to V$ mapping any
element  $ \sum_i p_i(h)
L_i\in A$  into
$  \sum_i p_i(\varepsilon) L_i \in V$
where $p_i(h)\in \Z[h]$ and $\{L_i\}_i$ are oriented links.
For $n\geq 0$, let $V^{\varepsilon}_n$
 be the image of $A_n$ under this homomorphism.
Clearly, $V^{\varepsilon}_n$ is the subgroup of $V$
generated by the resolutions of singular links with $n$ double points
where   we use the resolution, $r^{\varepsilon}$,
obtained from the one in Figure \ref {qres} by setting $h=\varepsilon$.
Corollary \ref {cor1} allows us
compute the   quotients associated with the filtration $V=
V^{\varepsilon}_0\supset V^{\varepsilon}_1\supset V^{\varepsilon}_2\supset
\ldots$. Namely,
  $V^{\varepsilon}/V^{\varepsilon}_{n}$ is a free abelian group  of
rank $n(n+1)/2$ freely
generated by the elements $  r^{\varepsilon}(G^l_{l+m}) \mod
V^{\varepsilon}_{n}$
with $ l + m < n$.
As above, for any oriented link $L$,  we have an expansion
\[L = \sum_{l,m\geq 0}  p^{\varepsilon}_{l,m} (L)\,
r^{\varepsilon}(G^l_{l+m})
 \in \varprojlim_n V/V^{\varepsilon}_n\]
with $p^{\varepsilon}_{l,m}(L)\in \Z$. Theorem \ref {deri} implies
that $$p^{\varepsilon}_{l,m} (L)= \nabla_{l,m}
(\varepsilon)=\varepsilon^{-l} (m!)^{-1}
(h^{l}\nabla_{l,0})^{(m)}(\varepsilon)$$ where $\nabla_{l,m}
=\nabla_{l,m} (L)\in \Lambda$. Thus, the numbers $
p^{\varepsilon}_{l,m}(L) $ are the coefficients in the expansion
of $\ h^{l} \nabla_{l,0} $ as
 a formal power series in $h-\varepsilon$:
$$\nabla_{l,0}(L)= \varepsilon^{l}h^{-l}\sum_{m\geq 0}
p^{\varepsilon}_{l,m}(L) (h-\varepsilon)^m.$$ By Theorem \ref
{thhomfly}, the numbers $\{p^{\varepsilon}_{l,m}(L)\}$ determine
the Homfly polynomial and are determined by it.

 3. In the sequel to this paper the authors
will consider higher Homfly skein modules and higher Kauffman
skein modules of 3-manifolds. See   also \cite{Josef2} and
\cite{Frohman} where   similar definitions of   higher skein
modules were suggested.

4. Assume we have a pair of oriented links $L_1$ and $L_2$ with
the same Homfly polynomial. Then by Theorem \ref{thhomfly} and
Theorem \ref{deri} we see that $\nabla_{l,m}(L_1) =
\nabla_{l,m}(L_2)$ for all $l,m$. So for the descending filtration
$\Lambda\otimes_{\Z[h]} A = \Lambda\otimes_{\Z[h]} A_0 \supset
\Lambda\otimes_{\Z[h]} A_1 \supset \ldots$ we have by Corollary
\ref{cor1} that
\[L_1 - L_2 \in \bigcap_n \left( \Lambda\otimes_{\Z[h]} A_n\right) .\]
It is well know that if an oriented link $L_2$ is the mutant of
some other oriented link $L_1$, then $\tilde{\nabla}(L_1) =
\tilde{\nabla}(L_2)$ (see \cite{Lick}). Hence it is certainly
clear from this that
\[\bigcap_n \left( \Lambda\otimes_{\Z [h]} A_n\right) \neq
\{0\}.\]

The paper is organised as follows. In Section 2 we prove that
$A_n/A_{n+1}$ is generated by the $h$-torsion and the generators
specified in Theorem \ref{Mainth}. In Section 3 we
introduce a certain quotient of the vector space  generated by chord diagrams
modulo the 4T-relation.
This quotient is used in Section 4 where we complete the proof of
Theorem \ref{Mainth}. In Section 5 we prove Theorem \ref {deri}.
In Section 6 we prove Theorem \ref  {thhomfly}.

\section{ The 8T-relation and generators of $A_n/A_{n+1}$}\label{gg}

We begin with a fundamental relation for singular links, which we
  call the 8T-relation.

\begin{figure}[htbp]
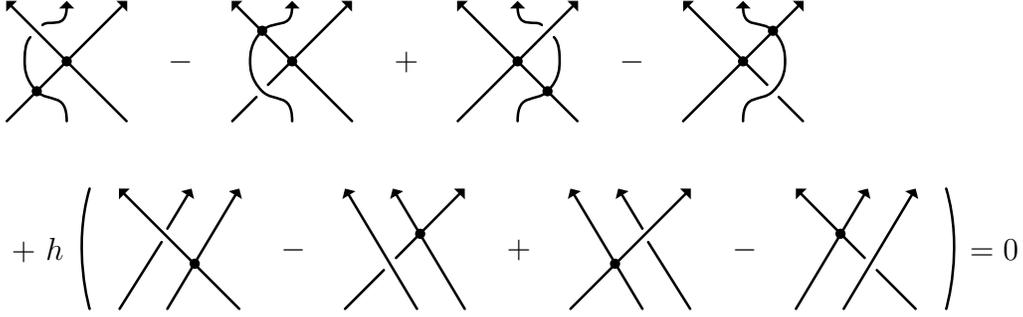

  \centertexdraw{\drawdim mm
\linewd 0.3 \arrowheadtype t:F \arrowheadsize l:1 w:2 \move(15
0)\xdota{0.8}\rmove(30 0)\xdotb{0.8}\rmove(30
0)\xdotc{0.8}\rmove(30 0)\xdotd{0.8} \move(30 8)\textref h:C v:C
\htext{\Large $-$} \move(60 8)\textref h:C v:C \htext{\Large $+$}
\move(90 8)\textref h:C v:C \htext{\Large $-$} \move(18 -25)
\vparent{0.8}\move(30 -25)\xdote{0.8}\rmove(30
0)\xdotf{0.8}\rmove(30 0)\xdotg{0.8}\rmove(30 0)
\xdoth{0.8}\rmove(12 0) \hparent{0.8} \move(14.5 -17)\textref h:R
v:C \htext{\Large $+\ h$} \move(45 -17)\textref h:C v:C
\htext{\Large $-$} \move(75 -17)\textref h:C v:C \htext{\Large
$+$} \move(105 -17)\textref h:C v:C \htext{\Large $-$} \move(135
-17)\textref h:L v:C \htext{\Large $= 0$} \move(142 -17) }
  \caption{The 8T-relation.}\label{R8T}
\end{figure}

\begin{prop}\label{8T}
We have the identity in Figure \ref{R8T}, once all double points
are resolved as in Figure \ref{qres}.
\end{prop}

It is understood that all eight local pictures in Figure \ref{R8T}
are completed by one and the same singular tangle to form eight
 singular links in the 3-sphere.
Alternatively, one may view the identity in Figure \ref{R8T} as a
formal relation between singular tangles which lies in the kernel
of the resolution map $r$.

\proof

Consider the strand leading from the second input to the second
output in the first four pictures. This strand contains one double
point and one over/under-crossing. Resolve this double point in
each of these four pictures. This yields an algebraic sum of eight
terms with coefficient $h^0$ and four terms with coefficient $h$.
The sum of eight terms vanishes while the sum of four terms is
exactly the opposite of the sum in the second row in Figure
\ref{R8T}.

\eproof

We shall need the following lemma

\begin{lem}\label {mult}
Let $a\in A_l$. Then disjoint union with $ha$ maps $A_n$ to
$A_{n+l+1}$.
\end{lem}

\begin{figure}[htbp]
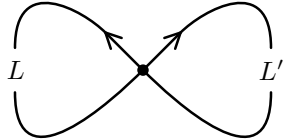

  \centertexdraw{
\drawdim mm \linewd 0.3 \arrowheadtype t:V \arrowheadsize l:2 w:2
\avec(5 5) \clvec(8 8)(17 13)(17 3) \move(17 -3) \clvec(17 -16)(4
-4)(0 0) \move(17 0)\textref h:C v:C \htext{$L'$} \move(-17
0)\textref h:C v:C \htext{$L$} \move(0 0) \avec(-5 5) \clvec(-8
8)(-17 13)(-17 3) \move(-17 -3) \clvec(-17 -16)(-4 -4)(0 0)
\move(0 12) \move(0 -12) \move(0 0) \fcir f:0 r:0.8 }
  \caption{The singular link N.}\label{nl1}
\end{figure}

\proof It suffices to prove that for any singular link $L$ with
$l$ double points and any singular link $L'$ with $n$ double
points, $h r(L\amalg L')\in A_{n+l+1}$. Consider the singular link
$N$ with $n+l+1$ double points obtained from $L$ and $L'$ as in
Figure \ref{nl1}. The result of resolving the double point in the
center is $-h L\amalg L'$, hence $h r(L\amalg L')=-r(N)\in
A_{n+l+1}$. \eproof

\begin{prop}\label{propgg}
For any $n\geq 0$,
\[A_n/A_{n+1} =  h{\text{-torsion}}\oplus\spanm_{\Z[h]}\{r(G^0_n), \ldots,
r(G^n_n)\}.\]

\end{prop}

\proof We first derive some consequences of the 8T-relation. For
any $(3,3)$-tangle with $n-1$ double points we complete the eight
local pictures in Figure \ref{R8T} with that tangle so as to
obtain eight singular links.
 The first four pictures in Figure \ref{R8T}
yield after resolution of double points elements of $A_{n+1}$,
which shall be ignored in the following calculations proceeding in
$A_n/A_{n+1}$. Thus we can complete the second row in Figure \ref{R8T} by
any $(3,3)$-tangle with $n-1$ double points and obtain a 4-term
relation in $A_n/A_{n+1}$.

 Let us now connect   the middle top strand to the bottom
left strand and add a
negative crossing at the bottom in the four pictures in the second row
of Figure \ref{R8T}. By the argument above, we obtain    a valid identity
in $A_n/A_{n+1}$, see Figure
\ref{PBR}.

\begin{figure}[htbp]
  \centertexdraw{\drawdim mm
\linewd 0.3 \arrowheadtype t:F \arrowheadsize l:1 w:2
\vparent{1}\rmove(15 0)\xdoteee{0.8}\rmove(30
0)\xdotfff{0.8}\rmove(30 0)\xdotggg{0.8}\rmove(30
0)\xdothhh{0.8}\rmove(10 0)\hparent{1} \move(-4 10)\textref h:R
v:C \htext{\Large $h$} \move(27 10)\textref h:C v:C \htext{\Large
$-$} \move(57 10)\textref h:C v:C \htext{\Large $+$} \move(86
10)\textref h:C v:C \htext{\Large $-$} \move(119 10)\textref h:L
v:C \htext{\Large $= 0 \mod A_{n+1}$} \move(137 0)}
  \caption{ \mbox{ }}\label{PBR}
\end{figure}

Observe that the first term in the equation in Figure \ref{PBR} is
in $A_{n+1}$ by Lemma \ref{mult}. Hence we obtain the {\em basic
relation} in $A_n/A_{n+1}$, see Figure \ref{BR}.

\begin{figure}[htbp]
  \centertexdraw{\drawdim mm
\linewd 0.3 \arrowheadtype t:V \arrowheadsize l:1 w:2 \textref h:c
v:c \xdot{0.8} \rmove(35 0)\vparent{0.8}\rmove(10
0)\xminusloop{0.8}\rmove(30 0)\xnulcirc{0.8}\rmove(10
0)\hparent{0.8} \move(-12 8)\textref h:R v:C \htext{\Large $h$}
\move(14 8)\textref h:L v:C \htext{\Large$=$} \move(32 8)\textref
h:R v:C \htext{\Large$h$} \move(57 8)\textref h:C v:C
\htext{\Large$+$} \move(87 8)\textref h:L v:C \htext{\Large $\mod
A_{n+1}$} \move(117 0)}
  \caption{The basic relation in $A_n/A_{n+1}$.}\label{BR}
\end{figure}

\begin{figure}[htbp]
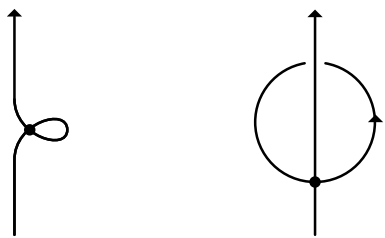

  \centertexdraw{
\drawdim mm \linewd 0.3 \arrowheadtype t:F \arrowheadsize l:1 w:2

\move(0 0) \rlvec(0 10) \clvec(0 15)(7 17)(7 14) \clvec(7 11)(0
13)(0 18) \ravec(0 12) \move(2 14) \fcir f:0 r:0.8

\move(40 0) \ravec(0 30) \move(40 15) \larc r:8 sd:-90 ed:80 \larc
r:8 sd:100 ed:270 \move(40 7) \fcir f:0 r:0.8 \move(48 15)
\ravec(0 1)
 }
  \caption{Singular tangles $T_1$ and $T_2$.}\label{sdp}
\end{figure}

To prove the proposition it is enough to show that for any
singular link $L$ with $n$ double points,
\[ h r(L) \,\in \,h \,\spanm_{\Z[h]}\{r(G^0_n),
\ldots, r(G^n_n)\}\mod A_{n+1}.\]
The basic relation implies that $hr(L) = \sum_{i=1}^{m}
hr(L_i)\mod A_{n+1}$ where $m=2^n$ and $L_1,\ldots,L_m$ are
singular links with $n$ double points such that all their
  double
points are as on the right-hand side of the equality in Figure \ref{sdp}. In
other
words, each $L_i$ is obtained from a non-singular link  by
inserting a certain number, say $l_i$, curls and attaching $n-l_i$
small unknotted circles meeting $L_i$ in one point as in Figure
\ref{sdp}.

To compute $ hr(L_i) \mod A_{n+1}$ we can use the same method as
in the usual recursive computation of the Conway polynomial of a
link. The role of the skein relation is played here by the fact
that we are computing
 modulo $A_{n+1}$. This shows that
 each $hr(L_i)\mod A_{n+1}$ expands as a linear combination over $\Z[h]$
 of certain $hr(L_i^j)$ where each $L_i^j$ is a disjoint union of singular
links of the form $G^s_t$. (We note that the order of any two
double points along the strand in Figure \ref{Gln} can be changed
in an arbitrary way without changing the resolution $r$. Indeed,
since any curl resolves to a disjoint union of an unknot times
$-h$ and the rest, we can move this unknot anywhere and reattach
it). By Lemma 2.2, if $L_i^j$ is disconnected then $hr(L_i^j)=0
\mod A_{n+1}$ so that we need to consider only connected $L_i^j$.
Then by the remarks above, $L_i^j=G_n^{l_i}$. Hence,
  $hr(L_i)= h p_i(h) r(G^{l_i}_n)\mod A_{n+1}$
where $p_i(h)\in \Z[h]$.
   This completes the proof of the
proposition.
 \eproof

 \rem \label{remConway} The arguments given in the proof of
  Proposition \ref{propgg} allow us to compute   the coefficients
  in the  expansion
  $$hr(L)= h\sum_{l=0}^{n} q_l  \,r(G^{l}_n)\mod A_{n+1}$$
  where $L$ is a singular link with $n$ double points and
 $q_l \in \Z[h]$.
 Denote by $\sing{L}$ the set of double points of $L$.
 For each   subset $X\subset \sing{L}$ denote by $L_X$ the non-singular link
obtained
from $L$ as follows:
 all double points of $L$ belonging to $X$ are replaced with   negative
crossings and all other double points of $L$ are smoothed. Then
  $$q_l=\sum_{\stackrel{{\scriptstyle X\subset \sing{L}}} {\card{X}=l}} \nabla (L_X)$$
  where $\nabla$ is the Conway polynomial of links.




\section{The 4TS-relation for chord diagrams}\label{s-uf}\label{Smoothforg}

\subsection{Chord diagrams.}
By a chord diagram we mean a finite family of oriented circles
with finitely many disjoint chords attached to them. Here a chord
connects either two distinct points of the same circle or two
points belonging to different circles.
 In our pictures, chords are
represented by fat dots, while the ordinary intersections of
strands should be ignored.

\begin{figure}[htbp]
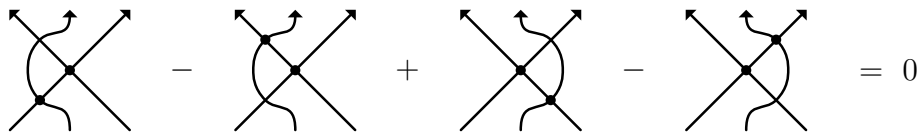

  \centertexdraw{\drawdim mm
\linewd 0.3 \arrowheadtype t:F \arrowheadsize l:1 w:2 \move(0
0)\xdotaa{0.8}\rmove(30 0)\xdotbb{0.8}\rmove(30
0)\xdotcc{0.8}\rmove(30 0)\xdotdd{0.8} \move(15 8)\textref h:C v:C
\htext{\Large $-$} \move(45 8)\textref h:C v:C \htext{\Large $+$}
\move(75 8)\textref h:C v:C \htext{\Large $-$} \move(105
8)\textref h:L v:C \htext{\Large $= \ 0$} \move(114 0)}
  \caption{The 4T-relation.}\label{4T}
\end{figure}

\begin{figure}[htbp]
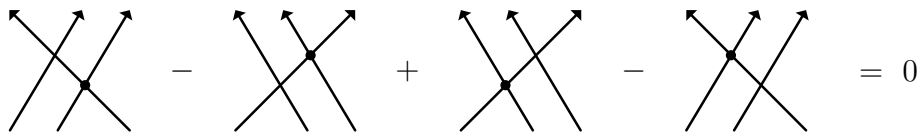

  \centertexdraw{\drawdim mm
\linewd 0.3 \arrowheadtype t:F \arrowheadsize l:1 w:2 \textref h:c
v:c \xdotee{0.8}\rmove(30 0)\xdotff{0.8}\rmove(30
0)\xdotgg{0.8}\rmove(30 0) \xdothh{0.8}\rmove(12 0) \move(15
8)\textref h:C v:C \htext{\Large $-$} \move(45 8)\textref h:C v:C
\htext{\Large $+$} \move(75 8)\textref h:C v:C \htext{\Large $-$}
\move(105 8)\textref h:L v:C \htext{\Large $= \ 0$} \move(114 0)}
  \caption{The 4TS-relation.}\label{4TS}
\end{figure}

For $n\geq 0$, let $ch_n$ be the  vector space over $\Q$ generated
by chord diagrams with $n$ chords.
We shall consider the following two quotients, $C_n$ and $D_n$, of $ch_n$:
$$C_n=ch_n/4TS\,\,\,\,\,\, {\text {and}}\,\,\,\,\,\,D_n=C_n/4T=ch_n/4T,4TS$$
where 4T is the standard 4-term relation for chord diagrams
shown in Figure \ref{4T} and 4TS is the   4-term relation for chord diagrams
shown in Figure \ref{4TS}.
The abbreviation 4TS should become more clear in Section
\ref{sf}.

The vector spaces $ch_n, C_n$, and $D_n$ have
  the structure of a module over the
 polynomial ring $\Q[u]$ on one variable $u$. The variable $u$ acts on a
chord diagram
$d$ by adding an oriented circle $g_0$ without chords:
 $ud=d\coprod g_0$.

Every $\Q[u]$-module $M$ has a completion $\hat M$ defined as the
limit of the projective system $M \leftarrow M/uM \leftarrow
M/u^2M \leftarrow...$. Note that $\hat M= \varprojlim_{e} M/u^e M$
is a module over the ring of formal power series $\Q[[u]]$. The
next proposition computes $\hat C_n$ and $\hat D_n$.

\begin{thm}\label{thmee}

 For any $n\geq 0$, the projection $C_n\to D_n$ induces an ismorphism
   $\hat C_n\to \hat D_n $. The completion
 $\hat C_n= \hat D_n $ is a free  $\Q [[u]]$-module
 of rank $n+1$ freely generated by the classes of the chord diagrams
 $g^l_n$, where $g^l_n$ is the underlying chord diagram of the singular link
$G^l_n$ shown in Figure \ref{Gln}, $l=0,1,\ldots,n$.

\end{thm}

The proof of this theorem given at the end of this section is based on a study
of three operators
 acting on chord
diagrams, specifically the operator adding  an isolated chord
and the   smoothing and forgetting operators.

 \subsection{Adding an isolated chord.}
 Let $T$ denote the operation on chord diagrams, which adds an
isolated chord. This operation is of course not well-defined on
the chord diagrams themselves, but is well-defined
provided we consider the diagrams modulo 4TS.

\begin{lem}\label{opT}
The operation of adding an isolated chord induces a well-defined
$\Q[u]$-linear homomorphism $T : C_{n}\to C_{n+1}$.
\end{lem}

\proof Consider the corollary of the 4TS-relation in Figure
\ref{4TS}, obtained by connecting the left top end to the left
bottom end on all four pictures. The first and second term cancel
and one obtains the equality of the middle two terms. This exactly
shows that an isolated chord can be moved from one place to any
other modulo 4TS. \eproof

\subsection{Smoothing and forgetting operators.}
\label{sf}  
\begin{figure}[htbp]
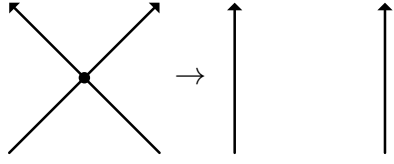

 \centertexdraw{\drawdim mm
\linewd 0.3 \arrowheadtype t:F \arrowheadsize l:1 w:2
\xdot{1}\move(30 0) \xnul{1} \move(14 10)\textref h:C v:C
\htext{\Large$\rightarrow$} }
  \caption{The smoothing of a chord.}\label{chsmooth}
\end{figure}

For each $n=1,2,...$, we define two homomorphisms $D_{n}\to
D_{n-1}$, called the smoothing and forgetting operators. The
smoothing operator, $S$,
 maps a chord diagram $d\in D_n$ into the  sum
$S(d)=\sum_c S(d;c)\in D_{n-1}$ where $c $ runs over all chords of
$d$ and $S(d; c)$ is $d$ with chord $c$ smoothed as shown in
Figure \ref{chsmooth}.
This operator maps
   4T and 4TS into 4TS and 0, respectively,
and defines therefore a $\Q[u]$-linear homomorphism $D_{n}\to
D_{n-1}$. Note that the 4TS-relation is nothing else than the
``4T-relation smoothed'', hence the name 4TS.

  The forgetting operator, $F$,    maps a chord diagram $d\in D_n$ into the
sum
$F(d)=\sum_c F(d;c)\in D_{n-1}$ where $c $ runs over all chords of
$d$ and $F(d; c)$ is $d$ with chord $c$ forgotten. It is easy to
check that $F$ maps both 4T
 and 4TS into $0$
and defines a $\Q[u]$-linear homomorphism $D_{n}\to D_{n-1}$.

The operator $T : C_{n}\to C_{n+1}$ constructed above induces
an operator $D_{n}\to D_{n+1}$ denoted by the same symbol $T$.
 We have the
following
commutation relations:
$$ SF-FS=0,\,\,\,\, ST-TS=u, \,\,\,\, FT-TF={\text{id}}.$$
This implies the   useful commutation relation
$ (S-uF)T=T(S-uF)$.

\subsection{Proof of Theorem \ref{thmee}.}
Both relations
4T
and 4TS are void for chord diagrams without chords. This implies
that $C_0=D_0$ is the free $\Q[u]$-module of rank 1 generated by
$g^0_0$. Clearly, $\hat C_0=\hat D_0=\Q[[u]] g^0_0$.

Assume  now that the claim of the theorem holds for   $n-1$
and prove it for $n$.
The proof consists of two parts. First
we show that the classes of the chord diagrams
$g^0_n,\ldots,g^n_n$ generate $\hat C_n$. Then we
show that
they are linearly independent in $\hat D_n$. These facts and the surjectivity
of the projection $\hat C_n\to \hat D_n$ would give the inductive step.

Let us prove that    $g^0_n,\ldots,g^n_n$ generate  $\hat C_n$.
By the inductive assumption, $\hat C_{n-1}$ is a free $\Q[[u]]$-module
freely generated by $g^0_{n-1}, \ldots ,g^{n-1}_{n-1}$.
Clearly, $T(g^l_{n-1})=g^{l+1}_n$.
Therefore the image of the   operator $\hat T:\hat C_{n-1}\to \hat C_n$
induced by $T:C_{n-1}\to C_n$
is  generated by $g^1_n,\ldots,g^n_n$. It remains to prove that
$g^0_n$ generates the $\Q[[u]]$-module $\hat C_n/ \hat T(\hat C_{n-1})$.

We shall show that any chord diagram $d$ with $n$ chords can be
expanded
 modulo 4TS and modulo $T(C_{n-1})$  in the form $d=\nu g^0_n +ud'$ where $\nu \in
\Z$ and $d'\in C_n$.
Iterating this expansion we obtain that $g^0_n$ generates the
$(\Q[u]/u^e)$-module
$$
C_n/(u^eC_n+T(C_{n-1}))=\hat C_n/(u^e\hat C_n+ \hat T(\hat C_{n-1})) $$
for any $e\geq 0$. This implies that
   $g^0_n$ generates $\hat C_n/\hat T(\hat C_{n-1})$ over $\Q[[u]]$.

The 4TS-relation implies the relation in Figure \ref{PBR} where we
remove $h$, ignore the over/under-crossings and interpret the four
terms as local pictures of chord diagrams. Therefore the basic
relation in Figure \ref{BR} with $h$ removed (again ignoring the
over/under-crossings) holds in $C_n/uC_n$. Quotienting by
$T(C_{n-1})$ we obtain the equalities in Figure \ref{cBR}.
Applying these equalities we can expand $d$ as a finite sum $\nu
g^0_n + \sum_j d_j$ where $\nu \in \Z$ and each $d_j\in C_n$ is a
disjoint union of several chord diagrams of type $g^0_m$ with
$m<n$. Each such $d_j$ belongs to $uC_n+T( C_{n-1})$. Indeed,
   the relation in Figure
\ref{cBR} implies that $$g^0_{m}\amalg g^0_{m'} = g^0_{m-1}\amalg
g^0_{m'+1}= \ldots = g^0_{0}\amalg g^0_{m+m'}=u g^0_{m+m'} \mod
T(C_{n-1}), 4TS$$ for all $m,m'\geq 0$. Hence, $d = \nu g^0_n +ud'
\mod T( C_{n-1})$ with $d'\in C_n$.

\begin{figure}[htbp]
  \centertexdraw{\drawdim mm
\linewd 0.3 \arrowheadtype t:F \arrowheadsize l:1 w:2 \xdot{0.8}
\rmove(30 0)\xnulcircv{0.8}\rmove(30 0)\xnulcirch{0.8} \move(13
7.5)\textref h:C v:C \htext{\Large $=$} \move(45 7.5)\textref h:C
v:C \htext{\Large $=$} \move(73 7.5)\textref h:L v:C \htext{\Large
$\mod T(C_{n-1}),4TS$} \move(110 0) }
  \caption{}\label{cBR}
\end{figure}

It remains to show that the classes of the chord diagrams
$g^0_n,\ldots,g^n_n$   are linearly independent in $\hat D_n$.
It follows from definitions that
\begin{equation}
  F(g^l_n)= l g^{l-1}_{n-1}+ u(n-l) g^{l}_{n-1}  \label{usui}
 \end{equation}
 for $l\geq 1$, and $F(g^0_n)=un g^{0}_{n-1}$.
Similarly,
$$S(g^l_n)= ul g^{l-1}_{n-1}+  (n-l) g^{l}_{n-1}  $$ for $l\geq 1$, and
$S(g^0_n)= n g^{0}_{n-1}$.
Now,
the $\Q[u]$-linear homomorphisms $S,F:D_n\to D_{n-1}$
induce $\Q[[u]]$-linear homomorphisms $\hat S, \hat F: \hat D_n\to \hat
D_{n-1}$. If there is a linear relation
$\sum_{l=0}^n k_l g^l_n=0$ in  $\hat D_n$ (where $k_l\in \Q[[u]]$) then applying
$\hat F$ and $\hat S$ we obtain two linear relations between the
classes of
$g^0_{n-1}, \ldots ,g^{n-1}_{n-1}$ in $\hat D_{n-1}$.
By the inductive assumption, these classes are linearly independent.
This gives two systems of linear equations on $\{k_l\}$:
first, $k_{l+1} (l+1)+k_l u (n-l) =0$
and second, $k_{l+1} u (l+1)+k_l (n-l) =0$ for $l=0,1,\ldots, n-1$.
The only solution is $k_l=0$ for all $l$. Thus,
the classes of
$g^0_n,\ldots,g^n_n$   are linearly independent in $\hat D_n$ which completes
the inductive step and the proof of the theorem.

 \eproof

\section{Proof of Theorem \ref{Mainth}}\label{linindep}

\subsection{Algebraic preliminaries.}  Let $K$ be a commutative algebra
over the field of rational numbers $\Q$. (In the sequel, $K$ will
be the polynomial ring $\Q[u]$). For any $K$-module $M$, denote by
$M[[v]]$ the set of formal power series on the
 variable $v$ with coefficients in $M$. We provide $M[[v]]$ with the
structure
of a module
over the ring of formal power series $K[[v]]$ in the obvious way.
We have $M\subset M[[v]]$: an element $a\in M$ is identified with
the formal power series $a+ 0\cdot v+0 \cdot v^2+\ldots$.

For $K$-modules $M_0,M_1,\ldots $, the product $\prod_{k\geq 0}
M_k$ is the $K$-module consisting of the series $a_0+a_1+\ldots$
with $a_k\in M_k, \,k\geq 0$. The addition and multiplication by
elements of $K$ are
 defined coordinate-wise. Applying this construction to the
 $K[[v]]$-modules $M_0[[v]],M_1[[v]],\ldots $ we obtain a $K[[v]]$-module
$ \prod_{k\geq 0} M_k [[v]]$. It is easy to observe that $
\prod_{k\geq 0} M_k [[v]]=(\prod_{k\geq 0} M_k) [[v]]$.

 Let  $M_0,M_1,\ldots $ be   $K$-modules provided with
  $K$-linear homomorphisms
$\alpha:M_{k}\to M_{k-1}$ for all $k\geq 1$. These morphisms
extend by linearity to $K[[v]]$-linear homomorphisms $
M_{k}[[v]]\to M_{k-1}[[v]]$ also denoted $\alpha$. We define a
$K[[v]]$-linear endomorphism, $\alpha$, of $ \prod_{k\geq 0} M_k
[[v]]$ by $$\alpha(\sum_{k\geq 0} a_k)= \sum_{k\geq 1}
\alpha(a_k)$$ where $a_k\in M_k[[v]]$ and $\alpha(a_k)\in
M_{k-1}[[v]]$.
 Finally, we define a $K[[v]]$-linear endomorphism  $e^{v\alpha} $ of
$ \prod_{k\geq 0} M_k [[v]]$ by $$e^{v\alpha} (\sum_{k\geq 0}
a_k)= \sum_{m\geq 0} \frac {v^m} {m!} \,\alpha^m(\sum_{k\geq 0}
a_k)= \sum_{k\geq 0} \left (\sum_{m\geq 0} \frac {v^m} {m!}
\,\alpha^m (a_{k+m })\right ).$$ It is easy to check that
$e^{v\alpha}$ is well defined and invertible
  with inverse $e^{-v\alpha}$.

\subsection{Framed singular links.} Recall that a framed link
 is a link provided with a homotopy class of non-zero normal vector fields.
 Let $\mathcal L^f$ denote the set of
the isotopy
 classes of framed oriented links in the 3-sphere.
 Denote by  $A^f$ the $\Z[h]$-module freely generated by the set $\mathcal
L^f$.
We
 define a filtration in $A^f$ using framed singular links as follows.
A framed singular link is a singular link (as defined in the introduction)
provided with a homotopy class of non-zero normal vector fields.
 In a neighbourhood of a double point
as in Figure \ref {doublept}  the vector field  should be
 orthogonal to the plane of the picture and
directed towards the reader. Note that we can keep the framing
when resolving a double point of a framed singular link as in
Figure \ref{qres}. (Here and below we use the standard convention
for the framings of links presented by link diagrams: the framings
are
 orthogonal to the plane of the pictures and  are
directed towards the reader).
Using the
formula in Figure \ref{qres}, we resolve each framed singular link $L$
with $n$ double points into a formal sum $r(L)\in A^f$ of $3^n$
terms. Denote by $A^f_n$   the $\Z[h]$-submodule of $A^f$ generated by $r(L)$
where $L$ runs over all framed singular links with $n$ double points.
Clearly, $A^f=A^f_0\supset A^f_1\supset A^f_2\supset \ldots$.
Forgetting the framing, we obtain a   projection $A^f\to A$ mapping each
$A^f_n$
onto $A_n$.

\subsection{The smoothed Kontsevich invariant.}\label{SmoothKontsevich}
The Kontsevich invariant of framed links is a mapping $\mathcal
L^f\to \prod_{k\geq 0} ch_k/4T$ where $ch_k/4T$ is the vector
space over $\Q$ generated by chord diagrams with $k$ chords modulo
4T (see e.g. \cite{Ko} and \cite{Bn}). Quotienting further by 4TS
we obtain a mapping $\mathcal L^f \to \prod_{k\geq 0} D_k$ where
$D_k=ch_k/(4T,4TS)$ is the $\Q[u]$-module considered in Section
\ref{Smoothforg}. We extend the latter mapping to an additive
homomorphism $z:A^f\to \prod_{k\geq 0} D_k [[v]]$ such that for
any $a\in A^f$ and any $ g(h)\in \Z[h]$ we have
\begin{equation}
z(g(h)a)= g(e^{v/2}-e^{-v/2})\, z(a)\label{uchi}
\end{equation}
(i.e. each entry of $h$ is traded for the formal power series
$e^{v/2}-e^{-v/2}$).

In the next   lemmas we shall combine $z$ with
  the $\Q[u]$-linear smoothing
and forgetting operators $S,F:D_k\to D_{k-1}$ defined in Section
\ref{Smoothforg}.   It is
  convenient to use the following formulas for
the endomorphisms $e^{vS} $  and $e^{vF} $ of
$ \prod_{k\geq 0} D_k [[v]]$ induced by $S$ and $F$:
for a chord diagram $d$,
 \begin{equation}e^{vS}(d)=\sum_C v^{\card{C}} S(d;C)
\,\,\,\, {\text { and }}\,\,\,\,e^{vF}(d)=\sum_C v^{\card{C}}
F(d;C)\label{expsf}
 \end{equation}
 where $C $ runs over all subsets of the set of chords of $d$ and
 $S(d;
C)$ (resp. $F(d; C)$) is $d$ with all the chords $c\in C$ smoothed
(resp. forgotten).

\begin{lem}\label{lem5.1}

For every $n\geq 0$, $$(e^{vS} z)(A^f_n) \subset \prod_{k\geq n}
D_k [[v]].$$ Moreover, if $L$ is a framed singular link in $S^3$
with $n$ double points then $$p_n e^{vS} z (r(L))=d(L) \mod v$$
where $p_n$ is the projection $\prod_{k\geq n} D_k [[v]] \to D_n
[[v]]$, $r(L)$ is the element of $A^f_n$ represented by $L$ and
$d(L)\in D_n$ is represented by the underlying chord diagram of
$L$.

\end{lem}

\proof

We begin with the second claim of the lemma. It is clear that
computing modulo $v$ we obtain $p_n e^{vS} z(r(L))=p_n z(r(L)) \mod
v$.
 Since the mapping $z$ maps each coefficient $h$ in the resolution
 $r(L)$ into $e^{v/2}-e^{-v/2}=v+v^3/24+\ldots$,
 the terms of this resolution with non-trivial powers of $h$ contribute 0 to
  $  z(r(L)) \mod v$.
Hence, $ z(r(L)) \mod v$ is just the Kontsevich invariant of the
standard Vassiliev resolution of $L$. The $n$-th term $p_n z(r(L))
\mod v$ of this invariant is well known to be $d(L)$.

To prove the first claim of the lemma, recall that the Kontsevich invariant
$z$
can be applied to a framed
tangle,
  the chords being attached to
  the underlying  1-manifold of the tangle. Here the orientation of this
1-manifold and
  the order of its endpoints is remembered while its embedding into the
3-space
is forgotten.
  The smoothing operator $S$ and the exponential
$e^{vS}$ extend to chord diagrams based on tangles and their
formal linear combinations over $\Q [[v]]$ in the obvious way.

   Let $X_+,X_-$, and $X_0$ be the framed oriented tangles
drawn in Figure \ref{fot}.
Set
$$Z_{\bullet}= z(X_+)- z(X_-) - (e^{v/2}-e^{-v/2})\,
z(X_0).$$
 We
claim that $$ p_0 e^{vS} ( Z_{\bullet})=0$$ where $p_0$ is the projection to the
module of
chord
diagrams with 0 chords. (The projection $p_0$ annihilates all
chord diagrams with at least one chord).

By definition, $$z(X_+)=\sum_{m\geq 0} \frac {t^m(X)} {2^m m!}
\,\,\,\, {\text { and }}\,\,\,\,z(X_-)=\sum_{m\geq 0} \frac
{(-1)^mt^m(X)} {2^m m!} $$ where $X$ is the underlying 1-manifold
of $X_+$ and $X_-$ and $t^m$ attaches to $X$ exactly $m$ parallel
chords connecting two components of $X$. By definition, $I_0=z(X_0)$
is the chord diagram consisting of two vertical arcs and no
chords. Then
\begin{eqnarray*}
 p_0 e^{vS} (Z_{\bullet})
 &=&  2 p_0 e^{vS}\left(
   \sum_{{\stackrel{{\scriptstyle m\geq 0}}{m \, \text{odd}}}}
 \frac {t^m(X)} {2^m m!}\right) -
(e^{v/2}-e^{-v/2}) p_0 e^{vS} (I_0)\\
 & = & 2\sum_{\stackrel{{\scriptstyle m\geq 0}}{m \, \text{odd}}} \frac {(vS)^m} {m!}
 (\frac {t^m(X)} {2^m m!} ) -
(e^{v/2}-e^{-v/2}) I_0 \\
 & = & 2\sum_{{\stackrel{{\scriptstyle m\geq 0}}{m \,
\text{odd}}}} \frac {v^m} {2^m m!} \,I_0 - (e^{v/2}-e^{-v/2}) I_0
= 0. \end{eqnarray*} Here we used the obvious equality $S^m( t^m
(X))=m! I_0$ for any odd $m$.

Let us prove the first claim of the lemma. It suffices to prove
that for any framed singular link $L$ in $S^3$ with $n$ double
points,
 $   e^{vS} z (r(L))$ expands as a sum of chord diagrams
with $\geq n$ chords. Let $B_1,\ldots,B_n$ be small 3-balls
surrounding the double points of $L$ and let $B$ be their
complement in the 3-sphere. Resolving $L$ we obtain an algebraic
sum of $3^n$ framed links which coincide in $B$ and represent a
framed tangle $\tau \subset B$. The
  Kontsevich invariant $z$ of  these links can be computed in two steps:
first
  compute $z$ for $\tau$ and for  the tangles sitting in
   $B_i,\, i=1,\ldots,n$, then
glue the resulting chord diagrams based on tangles along their
common endpoints on $\partial B$. Thus $ z (r(L))$ may be obtained
by gluing $ z (\tau)$ and $n$ copies of
  $ Z_{\bullet} $.
 Therefore by   Formula \ref{expsf}, $ e^{vS}
z (r(L))$ may be obtained by gluing $ e^{vS} z (\tau)$ and $n$
copies of
  $e^{vS} ( Z_{\bullet}) $.
By the result
  above, each of these $n$ copies
expands as a formal sum of chord diagrams with $\geq 1$ chords.
Therefore
  $   e^{vS} z (r(L))$ expands as a sum of chord diagrams
with $\geq n$ chords.

 \eproof

\begin{lem}\label{lem5.2}

For every $n\geq 0$, $$ e^{v(S-uF)} z (A^f_n) \subset
\prod_{k=0}^n (uv)^{n-k} D_k [[v]] \times \prod_{k\geq n+1} D_k
[[v]].$$ Moreover, if $L$ is a framed singular link in $S^3$ with
$n$ double points then $$ e^{v(S-uF)} z (r(L))=e^{-uvF} d(L) \mod
  v \prod_{k=0}^n  (uv)^{n-k} D_k [[v]] \times \prod_{k\geq n+1} D_k [[v]]. $$

\end{lem}

\proof

Since the operators $S$ and $F$ commute,
$e^{v(S-uF)}=e^{vS}e^{-uvF}=e^{-uvF}e^{vS}$.
By Lemma
\ref{lem5.1}, if $a\in A^f_n$ then $ e^{vS} z (s)=\sum_{i\geq n}
a_i$ with $a_i\in D_i [[v]]$. Moreover, for $a=r(L)$, we have $a_n
\in d(L) +v D_n [[v]]$. It remains to observe that $e^{-uvF}$ maps
$D_i [[v]]$ into
  $\prod_{k=0}^i  (uv)^{i-k} D_k [[v]]$.

 \eproof

\subsection{Proof of Theorem \ref{Mainth}.} Let $E_k=D_k/T(C_{k-1})$ be the
quotient of
$D_k$ by the subspace generated by chord diagrams with isolated chords.
It follows from  Theorem \ref {thmee} that
$\hat E_k$ is the free $\Q[[u]]$-module generated
by $g_k=g^0_k$. This implies that for all $m\geq 0$,
\begin{equation}u^mE_k/u^{m+1}E_k=\Q\,u^mg_k.\label{quoti}
 \end{equation}

Denote by $J$   the composition
of the mapping $e^{v(S-uF)} z:A^f\to \prod_{k\geq 0} D_k [[v]]$
and the projection   $$proj: \prod_{k\geq 0} D_k [[v]]\to
\prod_{k\geq 0} E_k [[v]].$$ Note  that
if a framed link $L'$ is obtained from a framed link $L$ by
inserting a $+1$
 framing twist
then $z(L')=e^{T/2} z(L)$ and therefore
 $$ J (L')= proj (e^{v(S-uF)} e^{T/2} z (L))=proj (e^{T/2} e^{v(S-uF)} z (L))
 =J (L) .$$
 Here the second equality follows from the fact that
 $S-uF$
commutes with $T$, so that $e^{v(S-uF)}$ commutes with
$e^{T/2}$. Thus, $J$ is framing-independent
and induces an additive
homomorphism from $A$ to $ \prod_{k\geq 0} E_k [[v]]$. Denote this
homomorphism by $j$.

It follows from Lemma \ref{lem5.2} that for
every $n\geq 0$,
 $$j(A_n) \subset \prod_{k=0}^n (uv)^{n-k} E_k
[[v]] \times \prod_{k\geq n+1} E_k [[v]].$$
 Therefore $j$ induces an
additive homomorphism
 $$j_n: A_n/A_{n+1}\to \prod_{k=0}^n
(uv)^{n-k} E_k [[v]]/u(uv)^{n-k}E_k [[v]]$$
 $$= \bigoplus_{k=0}^n
(uv)^{n-k} E_k [[v]]/u(uv)^{n-k}E_k [[v]].$$
Formula
\ref{quoti} implies that for each $a\in A_n/A_{n+1}$,
 \begin{equation}j_n(a)=
\sum_{k=0}^n (uv)^{n-k} j^k_n(a) g_k,\label{expanc}
 \end{equation}
 where $j^k_n(a)$ is a
uniquely defined element of $\Q[[v]]$. This gives $n+1$ additive
homomorphisms
 $$j^0_n,j^1_n,\ldots,j^n_n:A_n/A_{n+1}\to \Q[[v]]$$
satisfying (\ref{uchi}) for any $a\in A_n/A_{n+1}$ and any $
g(h)\in \Z[h]$.

In light of Proposition \ref{propgg}, to finish the proof of
Theorem \ref{Mainth}, we just need to show that the elements
$a^l=r(G_n^l)$ of $A_n/A_{n+1}$ represented by the singular links
$G^l_n$
 with $l=0,\ldots,n$ are linearly independent over
$\Z[h]$. It suffices to show that the $(n+1)\times (n+1)$-matrix
$(j^k_n (a^l))_{k,l}$ over $\Q[[v]]$ is non-degenerate. To this
end it suffices to compute this matrix modulo $ v$ and to show
that the resulting matrix over $\Q$ is non-degenerate.

 By
Lemma \ref{lem5.2},
 $$(uv)^{n-k} j^k_n(a^l) \,g_k= \frac {
(-uvF)^{n-k} (g^l_n)} {(n-k) !} \mod \, (u,v) (uv)^{n-k} E_k
[[v]]$$
 where $(u,v)$
is the ideal of $\Q[u][[v]]$ generated by $u $ and $
v$ and   $g^l_n$ is the underlying chord diagram of $G^l_n$.
Therefore $$j^k_n(a^l) \, g_k = (-1)^{n-k} ((n-k) !)^{-1} {
F^{n-k}}(g^l_n) \mod \, (u,v) .$$ It follows from Formula \ref {usui} that
 ${ F^{n-k}}(g^l_n)= 0 \mod\, u$ if $l< n-k$ and
${ F^{n-k}}(g^l_n)= l! \,g^{0}_k= l! \,g_k\mod\, u$ if $l=n-k$.
 This gives $j^k_n(a^l)= 0 \mod\, v$ if $l< n-k$ and
$j^k_n(a^l)=(-1)^{l} \mod\, v$ if $l=n-k$. Therefore the matrix
$(j^k_n (a^l))_{k,l}$
  is non-degenerate which completes the proof of the theorem.

\rem \label{ws} An easy calculation shows that
\[e^{vS}e^{T/2} = e^{vu/2}e^{T/2}e^{vS}.\]
From this we observe that $p_0 e^{vS}z(L)$ mod $uD_0[[v]]$ is
independent of the framing of the link $L$. Since $D_0 \cong
\Q[u],$ we have that $D_0[[v]]/uD_0[[v]]\cong \Q[[v]]$. So for an
oriented link $L$ we can now define $\sigma(L)\in \Q[[v]]$ by
\[\sigma(L) = p_0 e^{vS}z(L) \mod uD_0[[v]].\]
By Lemma \ref{lem5.1} we see that
\[\sigma(X_+) - \sigma(X_-) = (e^{v/2}- e^{-v/2})\sigma(X_0),\]
for any triple $(X_+,X_-,X_0)$ as in Figure 1. From this we see
that
\[\nabla(L)(e^{v/2}- e^{-v/2}) =  \sigma(L)/\sigma(G^0_0)\]
for any oriented link $L$, where $G^0_0$ is an oriented unknot.
 Let us now compute $\sigma(G^0_0)$.
Consider the singular link $G_1^1$. By Lemma \ref{lem5.1} we have
that $p_0 e^{vS}z(r(G_1^1))=0$. Now
\[z(r(G_1^1)) = e^{T/2}z(G^0_0) - e^{-T/2}z(G^0_0) - (e^{v/2}- e^{-v/2}) z(G_0^0)^2\]
so
\[e^{vu/2}p_0e^{vS}z(G_0^0) - e^{-vu/2}p_0e^{vS}z(G_0^0) = (e^{v/2}- e^{-v/2})
(p_0e^{vS}z(G_0^0))^2 \]
 in $D_0[[v]]\cong \Q[u][[v]]$. Hence
\[p_0e^{vS}z(G_0^0) = \frac{e^{vu/2}- e^{-vu/2}}{e^{v/2}- e^{-v/2}}\]
and therefore
\[\sigma(G_0^0) = \frac{v}{e^{v/2}- e^{-v/2}}.\]

If we write
\[\sigma = \sum_{k=1}^{\infty}\sigma_k v^k\]
then we can easily describe the weight system $w_k : ch_k \ra \Q$
which composed with $z$ gives $\sigma_k$. It is given by
\[w_k = p_0 \frac{S^k}{k!} \mod uD_0.\]
It is clear that the operator $p_0 e^{S}$ on a given chord diagram
simply just smoothes all chords in the diagram. The result is a
power of $u$ equal to the number of resulting components. Hence
$w_k$ on a chord diagram with $k$ chords is 1 if and only if
smoothing all chords results in a connected diagram. Thus we see
that $w_k$ vanishes on $4T, 4TS, 1T, u ch_k$ and takes the value
$1$ on the diagram $g^k_k\in ch_k$ for all $k$. By Theorem
\ref{thmee} any weight system on $ch_k$ with these properties
equals $w_k$.


\section{Proof of Theorem \ref {deri}}\label{Dif0}

In this section we study differentiation of link invariants
and prove Theorem \ref {deri}. A part of our results apply to links in arbitrary
3-manifolds and to arbitrary resolutions of double points.

\subsection{Differentiation of link invariants}\label{Diffinvar}
Let $K$ be a commutative ring endowed with a differential, i.e., with
an additive homomorphism
$x\mapsto x':K\to K$ such that $(xy)'=x'y+xy'$ for any $x,y\in K$.
For an oriented 3-manifold $M$, denote by $\mathcal A=\mathcal A(M,K)$ the free
$K$-module
generated by the isotopy classes of oriented links in $M$.
There is a unique additive homomorphism
$d:\mathcal A\to \mathcal A$ satisfying the following two conditions:

(i) $d$ maps the generators of $\mathcal A$ represented by oriented links   into
$0$;

(ii) for any $k\in K$, $a\in \mathcal A$, we have   $d(ka)=k'a + k d(a)$.

  We can compute $d$ explicitly as follows: if $a=\sum_i k_i
L_i\in \mathcal A$ where $k_i\in K$ and $\{L_i\}_i$ are oriented links in $M$
then
\begin{equation}
d(a)= \sum_i k'_i L_i\,\in \mathcal A.\label{lopo}
 \end{equation}
  Note that $d(ua)=ud(a)$ where
$u$ acts on $\mathcal A$ as the disjoint
union with an unknot.

  The dual $K$-module $\mathcal A^*={\text {Hom}}_{K}(\mathcal A,K)$
can be identified with the module of $K$-valued isotopy
invariants of oriented links in $M$.   For any $P\in \mathcal A^*$, consider
  the
  $K$-linear   homomorphism $d^*(P):\mathcal A\to K$
sending the  generator of $\mathcal A$ represented by an oriented link $L$ into
$(P(L))'$.  We can explicitly compute $d^*(P)$ as follows: if
$a=\sum_i k_i  L_i\in \mathcal A$ as above then $$d^*(P)(a)= \sum_i k_i\,
(P(L_i))' .$$
It is
clear that  $d^*:\mathcal A^*\to \mathcal A^*$ is an additive homomorphism such
that  $d^*(k
P)=k' P+ k d^*(P)$ for any
  $k\in K, P\in \mathcal A^*$.
The following lemma yields the fundamental relation between
 $d^*$
  and $d$.

 \begin{lem}\label {der1}
For any
  $a\in \mathcal A$ and   $P\in \mathcal A^*$,
  $$(P(a))'= d^*(P) (a)+ P( d(a)).$$
\end{lem}

\proof

If $a=\sum_i k_i L_i $ as above, then $$(P(a))'=(P(\sum_i
k_i L_i))'=
 \sum_i (k_i P(L_i))'$$
 $$=\sum_i k_i (P(L_i))'+\sum_i k'_i  P(L_i)$$
$$ =\sum_i k_i (P(L_i))'+P(\sum_i k'_i L_i)=d^*(P) (a)+ P(
d(a)).$$ \eproof

We describe now the behaviour of any Vassiliev-type filtration in
$\mathcal A=\mathcal A(M,K)$ under the differential $d$. Fix a
finite formal linear combination $\sum_j k_j T_j$ where
 $k_j\in K$ and
each $T_j$ is a tangle in the 3-ball with two inputs and two
outputs. Consider a resolution $R$ of a double point (Figure \ref
{doublept}) defined by $R(X_\bullet)=\sum_j k_j T_j$. In this way
  we resolve each singular link $L\subset M$
  into a formal sum $R(L)\in \mathcal A$.
 Denote by $\mathcal A_n$   the $K$-submodule of $\mathcal A$ generated by
$R(L)$
where $L$ runs over all singular links with $n$ double points in $M$.
Clearly, $\mathcal A=\mathcal A_0\supset \mathcal A_1\supset   \ldots$.

\begin{lem}\label {der2}
For each $n\geq 0$, $d(\mathcal A_{n+1})\subset \mathcal A_n$.
\end{lem}

\proof
Observe first that the definition of the differential $d$ extends word
for word to linear combinations of oriented
tangles in oriented 3-manifolds with coefficients in
$K$ (use Formula \ref {lopo}).
Note that the usual gluing of tangles extends by linearity to their linear
combinations. It is clear that
if $EF$ is the result of gluing of two tangles
 (or linear combinations there of) $E,F$ then
$d(EF)=d(E)  F+E d(F)$.

To prove the lemma, it is enough to prove that for any singular link $L
 $ with
$n+1$ double points, $d(R(L)) $ is a linear combination of
 the resolutions of  singular links
 with $n$ double points.  We shall prove a more general claim: for any singular
tangle $L$ with $n+1$ double points, $d(R(L)) $ is a linear
combination of the resolutions of singular tangles with $n$ double
points. The proof goes by induction. For $n=0$ the claim is
obvious. Assume that the claim holds for $n< N$ and prove it for
$n=N$. Consider a singular tangle $L$ with $N+1$ double points.
Choose a double point $x$ of $L$ and  split $L$ into
two pieces: the singular tangle $X_{\bullet}$
  in a 3-ball neighborhood of $x$ and the
complementary singular tangle $\tau$ in the complement of this
3-ball. It follows from definitions that $R(L)=\sum_j k_j T_j
R(\tau)$. Note that $d(T_j)=0$. Hence $$ d(R(L))=\sum_j k'_j T_j
R(\tau)+ \sum_j k_j T_j d(R(\tau))$$ $$=\sum_j k'_j R(T_j \tau)+
R(X_{\bullet}) d(R(\tau)).$$ Since $\tau$ is a singular tangle
with $N$ double points, the inductive assumption implies that
$d(R(\tau))$ is a linear combination of the resolutions of
singular tangles with $N-1$ double points. Hence, $R(X_{\bullet})
d(R(\tau))$ is a linear combination of the resolutions of singular
tangles with $N$ double points. This proves the inductive step.

 \eproof

\begin{cor}\label {der2.5}
If $P\in \mathcal A^*$ annihilates $\mathcal A_{n}$ with $n \geq 0$ then
$d^*(P)$
annihilates $\mathcal A_{n+1}$.
\end{cor}

\proof For any $a\in \mathcal A_{n+1}\subset \mathcal A_n$, we have $P(a)=0$. By
the
previous lemma, $P(d(a))=0$. By Lemma \ref {der1},
 $d^*(P) (a)=(P(a))'- P( d(a))=0$.

\begin{cor}\label {der2.6}
For any
$n\geq 1$, the differential   $d:\mathcal A\to \mathcal A$ induces
an additive homomorphism
$\mathcal A/\mathcal A_{n+1}\to \mathcal A/\mathcal A_{n}$.
Its restriction
 $\mathcal A_n/\mathcal A_{n+1}\to \mathcal A_{n-1}/\mathcal A_{n}$
 is  $K$-linear.
\end{cor}

\subsection{Proof of Theorem \ref {deri}}
We apply the differentials $d$ and $d^*$ introduced above in the
case $M=S^3$ and $K=\Lambda=\Z[h,h^{-1}]$ with usual
differentiation of Laurent-polynomials in $h$. Note that in this
case $\mathcal A=\Lambda \otimes_{\Z[h]} A$ and $\mathcal
A/\mathcal A_n= \Lambda \otimes_{\Z[h]} (A/A_n)$ for all $n\geq
0$.

Recall that $$\P_{l,m}=(-h)^{l}\nabla_{l,m}\in {\text
{Hom}}_{\Z[h]}(A,\Lambda)= {\text {Hom}}_{\Lambda}(\mathcal
A,\Lambda)=\mathcal A^*.$$ We should prove that
\[\P_{l,m}= (m!)^{-1} (d^*)^m (\P_{l,0}). \]
It is enough to prove that $$d^* (\P_{l,m-1})= m \P_{l,m} $$ for
all $l\geq 0, m\geq 1$. By Corollary \ref{der2.5}, both sides
annihilate $\mathcal A_{l+m+1}$ and determine $\Lambda$-linear
homomorphisms $ \mathcal A/\mathcal A_{l+m+1} \to \Lambda$. It
suffices to verify that these homomorphisms coincide on the
generators $u^s r(G^0_t)$ of $\mathcal A/\mathcal A_{l+m+1}$ where
$s+t\leq l+m$. Set $a_t=r(G^0_t)\in \mathcal A$. By definition
(cf. Remark \ref{remIntrod} 1.), $ \P_{l,m}(u^s a_t)= \delta_l^s
\delta_m^t$. By Lemma \ref {der1}, $$d^* (\P_{l,m-1})(u^s
a_t)=(\P_{l,m-1} (u^s a_t))'- \P_{l,m-1} (u^s d(a_t))=-\P_{l,m-1}
(u^s d(a_t)).$$ If $t=0$, then $d(a_t)=d(r(G_0))=d(G_0)=0$ and
$$-\P_{l,m-1} (u^s d(a_t))=0= m\delta_l^s \delta_m^t. $$
  Assume that  $t\geq 1$.
  It is clear that $G_t$ is the closure
of the $t$-th power of the singular tangle $T_2$ drawn in Figure \ref
{sdp}. Therefore $r(G_t)$ is the closure of
$r(T_2^t)$. It follows from definitions  that $d(r(T_2))= - I$ where $I$ is
the
 unknotted vertical strand oriented upwards. Therefore,
$$ d(r(T_2^t))=d( (r(T_2))^t)=-t(r(T_2))^{t-1}=-t r(T_2^{t-1}).$$
Taking the closures, we obtain that $d(r(G_t))=-t\, r(G_{t-1}) $.
Thus, $$-\P_{l,m-1} (u^s d(a_t))= t \P_{l,m-1} (u^s a_{t-1})
=t\delta_l^s \delta_{m-1}^{t-1}=m\delta_l^s \delta_m^t.$$

\rem Weight systems can be obtained for the derivatives of the
Conway polynomial using the weight systems for the Conway
polynomial described in Remark \ref{ws}. This will be treated
elsewhere.

\section{Proof of Theorem \ref {thhomfly}}\label{Homomfly}


  For any oriented link $L$, we define a two variable
formal power series $P(L) \in \Z[h][[u]] $ by
\[P(L)(h,u)=\sum_{l=0}^\infty (-h)^{l}\nabla_{l,0}(L)(h) u^l.\]
This definition extends to singular links by  $P(L) = P(r(L))$.

The function $L\mapsto P(L)$ is $u$-linear, i.e., $P(u L) = uP(L)$
for any oriented link or singular link $L$. Furthermore,
$u^{n}|P(A_n)$, since $\nabla_{l,0} (A_n)=0$ for $l<n$. To prove
Theorem \ref {thhomfly} we establish the following theorem.

\begin{thm}\label{Homflyres}
There is a   pair of formal power series $\alpha, \beta\in
\Z[[u]]$ such that $P$ satisfies the skein relation
\[P(X_+) = (1+h\alpha)P(X_-) + (h+h\beta)P(X_0)  \]
for any three oriented links
$X_+, X_-,X_0$ coinciding outside a 3-ball
and looking
as in Figure \ref {fot} inside this ball.
The power series $\alpha$ and $\beta$ are described in Proposition
\ref{formula}.
\end{thm}

It follows from this theorem that $P$ is a reparametrised version of the
Homfly
link polynomial.

\begin{prop}\label{formula}
There is a unique pair of formal power series $\alpha, \beta\in
\Z[[u]]$ which satisfy the equations
 \begin{eqnarray}\alpha = -u +
  \beta(\alpha + \beta u)\label{uuu1}\\
  \beta = \alpha(\alpha + \beta u)\label{uuu2}
 \end{eqnarray}
and such that $\alpha = -u \mod u^2$ and $\beta = 0 \mod u^2$.
\end{prop}

\proof

Let us assume first that $\alpha$ and $\beta$ exist and show their
uniqueness. Since $\alpha=-u \mod u^2$, the formal power series
$\alpha+u\beta$ is divisible by $u$ so there exists $\gamma \in
\Z[[u]]$ such that $ \alpha+u\beta=u\gamma$. Clearly, the free
term of $\gamma$ is $-1$, so that $\gamma$ is invertible in
$\Z[[u]]$. Multiplying Formula \ref{uuu2} by $u$ and adding it to
Formula \ref{uuu1}, we obtain $u\gamma =-u+ u\gamma
(\beta+u\alpha)$ which is equivalent to $\gamma =-1+ \gamma
(\beta+u\alpha)$. This implies that $\beta+u\alpha=1+\gamma^{-1}$
so that $\beta=1+\gamma^{-1}-u\alpha$. Substituting this
expression for $\beta$ in the formula $ \alpha+u\beta=u\gamma$ we
obtain a linear equation on $\alpha$ which yields
 \begin{equation}
 \alpha= u(1-u^2)^{-1} (\gamma-\gamma^{-1}-1).\label{f1}
 \end{equation}
 We can also
determine $\beta$ from the equality $ \alpha+u\beta=u\gamma$. This
gives
\begin{equation}
 \beta= (1-u^2)^{-1} (1-u^2\gamma+\gamma^{-1}).\label{f2}
 \end{equation}

Substituting these expressions into \ref{uuu1} and \ref{uuu2}, we
easily observe that these two equations are equivalent to the
following equation on $\gamma$:
\begin{equation}
 u^2\gamma^3-(u^2+1) \gamma
-1=0.\label{f3}
 \end{equation}
 Now, expanding $\gamma=-1+\sum_{k\geq 1} a_k u^k$ with
$a_k\in \Z$ we inductively compute the coefficients of $\gamma$ from \ref
{f3}.
Hence there is only one formal power series $\gamma$ satisfying
this equation. (In fact, $\gamma$ is a formal power series in
 $u^2$). This proves uniqueness of $\alpha$ and $\beta$.
 Conversely, definining $\alpha,\beta, \gamma$ by Formulas
 \ref{f1} - \ref{f3}   we obtain $\alpha, \beta$ satisfying the conditions of
the proposition.

 \eproof

The key ingredient in the proof of   Theorem \ref{Homflyres}
is the following local relation.

\begin{lem}\label{Basicl}
Let $\alpha, \beta\in \Z[[u]]$ be as   in Proposition
\ref{formula}. Then
 \begin{equation}
 P(X_\bullet) = h \alpha P(X_-) + h\beta P(X_0) \label{Pdoublept}
 \end{equation}
 where $X_\bullet$ is the double point as in Figure \ref {doublept}.
\end{lem}

We note that   Theorem \ref {Homflyres} follows directly from this lemma,
since
by
definition
$$P(X_\bullet) = P(r(X_\bullet) )= P(X_+) - P(X_-) - h P(X_0).$$

\begin{figure}[htbp]
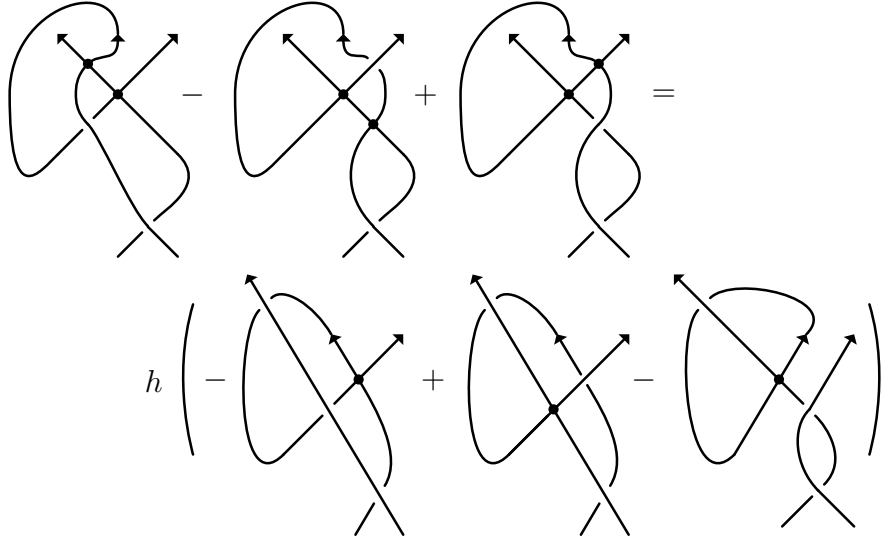

  \centertexdraw{\drawdim mm
\linewd 0.3 \arrowheadtype t:F \arrowheadsize l:1 w:2
\xdotbbb{0.8}\rmove(30 0)\xdotccc{0.8}\rmove(30 0)\xdotddd{0.8}
\move(10 8)\textref h:C v:C \htext{\Large $-$} \move(41 8)\textref
h:C v:C \htext{\Large $+$} \move(71 8)\textref h:L v:C
\htext{\Large $=$} \move(10 -40)\vparent{1}\rmove(20
0)\xdotfff{0.8}\rmove(30 0)\xdotggg{0.8}\rmove(30
0)\xdothhh{0.8}\rmove(10 0)\hparent{1} \move(6 -30)\textref h:R
v:C \htext{\Large $h$} \move(13 -30)\textref h:C v:C \htext{\Large
$-$} \move(42 -30)\textref h:C v:C \htext{\Large $+$} \move(70
-30)\textref h:C v:C \htext{\Large $-$} \move(0 26) }
  \caption{The 6T-relation.}\label{6T}
\end{figure}

\proof By connecting the middle top strand to the bottom left
strand and by adding a negative crossing at the bottom of the 8
terms in the 8T-relation (as we did in Section \ref{gg}), we
obtain, after canceling the first and the fifth term, the
``6T-relation'' in Figure \ref{6T}, which holds in $A$ once all
double points have been resolved.

We now claim that $P$ vanishes on any singular link obtained by
completing the singular $(1,1)$-tangle $T_2$ in Figure \ref{sdp}
by any singular (1,1)-tangle. Let $L$ be such a singular link with
$n$ double points. We apply the argument
  given in the proof of Proposition \ref{propgg}
  to all double points of $L$ except to the one of $T_2\subset L$.
  This and the   proof of
Lemma \ref{mult} show that $hr(L)\in A$ may be expanded as a finite sum
 $\sum_{l=0}^{n-1} h q_l r(G^l_n) + \sum_{j} h a_j r(L_j)$ where $q_l,a_j\in
\Z[h]$
 and each $L_j$ is a singular link  with $\geq n+1$ double points
  obtained as a
completion of $T_2$ by a singular (1,1)-tangle with $\geq n$ double points.
Repeating this argument inductively, we see that for any $N\geq n$
there is a finite
 expansion
$$hr(L)=
\sum_{m=n}^N\sum_{l=0}^{m-1} h q_{l,m}\, r(G^l_m) + \sum_{j} h a_{j,N}
\,r(L_{j,N})\,\in
A$$
 where $q_{l,m}, a_{j,N}\in \Z[h]$
and $L_{j,N}$ are singular links with $\geq N+1$ double points
  obtained by a
completion of $T_2$. Note that $P(G^l_m)=P(r(G^l_m))=0$ for $ l\leq m-1$
  and $P(L_{j,N})$ is divisible by $u^{N+1}$.
 Therefore   $P(L)$ is divisible by any
power  of $u$, hence it vanishes.

From the above we conclude that $P$ satisfies the relation in
Figure \ref{expXbullet}.
\begin{figure}[htbp]
  \centertexdraw{\drawdim mm
\linewd 0.3 \arrowheadtype t:F \arrowheadsize l:1 w:2
\vparent{0.8} \rmove(10 0) \xdot{0.8} \rmove(10 0) \hparent{0.8}
\rmove(30 0) \vparent{0.8} \rmove(10 0) \xminus{0.8} \rmove(10 0)
\hparent{0.8} \move(24 -25) \vparent{0.8} \rmove(11 0)\xp{0.8}
\rmove(11 0) \hparent{0.8} \rmove(14 0) \vparent{0.8}\rmove(11
0)\xpminus{0.8}\rmove(11 0) \hparent{0.8} \rmove(14 0)
\vparent{0.8}\rmove(11 0)\xpplus{0.8}\rmove(11 0) \hparent{0.8}
\move(16 -27) \vparent{1} \rmove(106 0) \hparent{1} \move(-4
8)\textref h:R v:C \htext{\Large $P$} \move(25 8)\textref h:L v:C
\htext{\Large $=$} \move(46 8)\textref h:R v:C \htext{\Large $-h u
P$} \move(12 -17)\textref h:R v:C \htext{\Large $+\ h^{-1}$}
\move(18.5 -17)\textref h:C v:C \htext{\Large $P$} \move(52.5
-17)\textref h:C v:C \htext{\Large $-\ P$} \move(88.5 -17)\textref
h:C v:C \htext{\Large $+\ P$} \move(109 0) }
  \caption{ \mbox{ }}\label{expXbullet}
\end{figure}
From this formula we see that $P(X_\bullet) = -hu P(X_-) \mod
u^2$, since the last three terms are divisible by $u^2$. This
agrees with (\ref{Pdoublept}) modulo $u^2$. We shall now prove
Formula (\ref{Pdoublept}) by proving it mod $u^{k}$ by induction
on $k$. Let us assume that (\ref{Pdoublept}) holds modulo $u^k$.
Since both $\alpha$ and $\beta$ are divisible by $u$, we can apply
our mod $u^k$ formula for $P(X_\bullet)$ to each of the double
points in the last three terms in Figure \ref{expXbullet}, and
obtain a formula for $P(X_\bullet)$ mod $u^{k+1}$. When we do
that, we obtain
\[P(X_\bullet) = h(-u + \beta( \alpha+\beta u ))P(X_-) + h\alpha(\alpha +
\beta u) P(X_0) \mod u^{k+1}.\] By Proposition \ref{formula},
 this is exactly Formula (\ref{Pdoublept})   modulo
$u^{k+1}$. \eproof


\end{document}